\documentclass[acmtoms]{acmtrans2m}

\usepackage{etex}

\makeatletter
\def\ps@myheadings{\let\@mkboth\@gobbletwo
\def\@oddhead{\hbox{}\hfill \small\sf \rightmark\hskip 19pt{\Large$\cdot$}\hskip 17pt\mypage}
\def\@oddfoot{\hbox{}\hfill\tiny\@runningfoot}
\def\@evenhead{\small\sf\mypage \hskip 17pt{\Large$\cdot$}\hskip 19pt\leftmark\hfill \hbox{}}
\def\@evenfoot{\tiny\@runningfoot\hfill\hbox{}}
\def\sectionmark##1{}\def\subsectionmark##1{}}
\def\@runningfoot{}
\def\runningfoot{\def\@runningfoot{}}
\def\@firstfoot{}
\def\firstfoot{\def\@firstfoot{}}
\def\ps@titlepage{\let\@mkboth\@gobbletwo
\def\@oddhead{}\def\@oddfoot{\hbox{}\hfill
\tiny\@firstfoot}\def\@evenhead{}\def\@evenfoot{\tiny\@firstfoot\hfill\hbox{}}}
\makeatother

\usepackage{numprint}
\usepackage{easybmat,easyeqn}
\usepackage[boxruled,linesnumbered]{algorithm2e}
\usepackage{graphicx}
\usepackage{pgf,pgfplots,tikz}
\usetikzlibrary{arrows,patterns,plotmarks}
\usepackage{bm,bbm,url}

\newcommand{\assign}{\leftarrow}

\newcommand{\zero}{\bm{0}}

\newcommand{\AAmat}{\bm{\mathcal{A}}}
\newcommand{\AAscal}{\mathcal{A}}

\newcommand{\MMscal}{\mathcal{M}}

\newcommand{\PPscal}{\mathcal{P}}

\newcommand{\ZZscal}{\mathcal{Z}}

\newcommand{\Xor}{\mathsf{Xor}}
\newcommand{\And}{\mathsf{And}}

\newcommand{\popC}{\mathtt{\,popC}}

\newcommand{\Amat}{\bm{A}}
\newcommand{\Bmat}{\bm{B}}
\newcommand{\Cmat}{\bm{C}}
\newcommand{\Dmat}{\bm{D}}
\newcommand{\Emat}{\bm{E}}

\newcommand{\Hmat}{\bm{H}}
\newcommand{\Jmat}{\bm{J}}

\newcommand{\Imat}{\bm{I}}
\newcommand{\Lmat}{\bm{L}}
\newcommand{\Mmat}{\bm{M}}

\newcommand{\Pmat}{\bm{P}}
\newcommand{\Qmat}{\bm{Q}}
\newcommand{\Rmat}{\bm{R}}

\newcommand{\Umat}{\bm{U}}

\newcommand{\Ymat}{\bm{Y}}
\newcommand{\Zmat}{\bm{Z}}

\newcommand{\cvec}{\bm{c}}
\newcommand{\dvec}{\bm{d}}
\newcommand{\evec}{\bm{e}}

\newcommand{\rvec}{\bm{r}}
\newcommand{\svec}{\bm{s}}

\newcommand{\FF}{\mathbbm{F}}
\newcommand{\NN}{\mathbbm{N}}

\newcommand{\AND}{\odot}
\newcommand{\XOR}{\oplus}
\newcommand{\rk}{\textsf{rank}}
\newcommand{\matSet}{\textsf{Mat}}

\newcommand{\bigO}{\mathcal{O}}

\newtheorem{theorem}{Theorem}[section]

\newtheorem{lemma}[theorem]{Lemma}

\newdef{remark}[theorem]{Remark}
\newtheorem{example}[theorem]{Example}

\npdecimalsign{.}

\begin{document}

\markboth{E. Bertolazzi, A. Rimoldi }{Fast matrix decomposition in $\FF_2$}

\title{Fast matrix decomposition in $\FF_2$}

\author{%
Enrico Bertolazzi\\
Department of Mechanical and Structural Engineering -- University of Trento, Italy, 
\\ \\
Anna Rimoldi\\
Department of Mathematics -- University of Trento, Italy.}

\begin{abstract}
In this work an efficient algorithm to perform a block decomposition for 
large dense rectangular matrices with entries in $\FF_2$ is presented.
Matrices are stored as column blocks of row major matrices 
in order to facilitate rows operation and matrix multiplications 
with block of columns.
One of the major bottlenecks of matrix decomposition is the
pivoting involving both rows and column exchanges.
Since row swaps are cheap and column swaps are order of magnitude
slower, the number of column swaps should be reduced as much as possible.
Here is presented an algorithm that completely avoids the column 
permutations.
An asymptotically fast algorithm is obtained by combining the four Russian algorithm
and the recursion with Strassen algorithm for matrix--matrix multiplication.
Moreover optimal parameters for the tuning of the algorithm 
are theoretically estimated and then experimentally verified.
A comparison with the state of the art public domain 
software SAGE shows that the proposed algorithm is generally faster.

\end{abstract}

\category{}{Numerical analysis}{Computations in finite fields}

\terms{Algorithms, Performance}


\maketitle




\section{Introduction}
An important tool in linear algebra is the matrix decomposition, which expresses
a (rectangular) matrix  as a product of two or more 
simpler matrices.
Such decompositions are used for easy computation of rank, null space, 
and solving linear system and related problem.

There are many well-known algorithms for matrix decomposition defined in any field, finite or not.
A common approach consists in reducing a matrix to the row-echelon 
form by row operations \cite{Meyer:2000}.
Once the row-echelon form is obtained, the rank will be equal to the number of non-zero 
rows and null space can be easily computed.
Gauss $LU$ decomposition \cite{Golub:1996} can be also used
to solve linear systems (when the matrix is square and full rank) 
or compute the rank.
Applied to rectangular or rank deficient matrices, 
it is costly as the computation of the row-echelon form.
In fact, Gauss decomposition of a matrix $\Amat$ produces an $LU$ factorization,
i.e. $\Pmat\Amat\Qmat=\Lmat\Umat$ where $\Pmat$ and $\Qmat$ are permutation matrices,
$\Lmat$ is the lower triangular matrix
and $\Umat$ is the row-echelon form of $\Amat$ (up to column permutations).
The asymptotic cost of a naive implementation of $LU$ decomposition for a dense 
$n \times n$ matrix is $\bigO(n^3)$.
However such a cost can be reduced using a combination of recursion
and matrix-matrix multiplication.
For example, using matrix-matrix multiplication (MMM)
in the construction of the factorization, the asymptotic cost
may be reduced by using fast MMM algorithms.
The problem of fast matrix-matrix multiplication is still
under development.

The \emph{naive} MMM algorithm is based on the classical definition
of the multiplication of two matrices; 
its cost is $n^3$ multiplications and $n^2(n-1)$ additions and so 
we classify the naive algorithm as an $O(n^3)$ algorithm.

\emph{Strassen} matrix-matrix multiplication algorithm \cite{Strassen:1969}
-- which asymptotic cost is $\bigO(n^{\log_27})$ in any field --
uses only seven scalar multiplications (instead of the usual eight) 
to multiply $2 \times 2$ matrices. In fact, as proved in \cite{Winograd:1971}, 
Strassen's algorithm is optimal for $2 \times 2$ matrices.
Further asymptotic improvements \cite{Coppersmith:1990} 
can be obtained to perform multiplication of larger matrices.
Hybrid algorithms incorporate Strassen and Winograd variants recursively 
to achieve high performance on large matrices
\cite{Huss-Lederman:1996,Higham:1990,Douglas:1994,Kaporin:1999}.
The asymptotic cost $\bigO(n^{\log_27})$ means that for a large enough $n$,
Strassen's algorithm should theoretically perform multiplication 
significantly faster than the naive algorithm.
However, asymptotic cost means that the actual cost of standard $LU$ decomposition is about 
$C_1 n^3$ while the actual cost of Strassen multiplication is about $C_2 n^{\log_27}$ 
where $C_2 \gg C_1$. Therefore, the use of Strassen algorithm is convenient only for large $n$.
Strassen algorithm is recursive so that normally the recursion is terminated 
when the cost of recursion is larger than the classical matrix-matrix multiplication.
That happens when $C_2 n^{\log_27} \approx C_1n^3$, i.e. $n\approx \exp(\ln(C_2/C_1)/(3-\log_27))$.
For instance, if $C_2/C_1\approx 10$ we have $n\approx 150000$ while
in case $C_2/C_1\approx 5$ we have $n\approx 4000$.
In practice, the computation of the switching point must consider additional 
costs and it is implementation dependent. For a detailed analysis see for example
\cite{Huss-Lederman:1996,Higham:1990}.

The efficient computation of MMM can be further improved in case of finite fields. 
In particular, in case of $\FF_2$, the 
\emph{Method of four Russian for Multiplication} (M4RM) is a fast MMM algorithm, 
which cost is $\bigO(n^3/\log n)$ \cite{Arlazarov:1970,aho:1974,Albrecht:2010}.
Its asymptotic cost is better than classical matrix-matrix multiplication;
but it is worse than recursive Strassen's algorithm.
However, if the actual cost of M4RM is about $C_3 n^3/\log n$,
we have that $C_3 \ll C_1$ so that is competitive for not too big matrices.
A combination of Strassen and M4RM is a good compromise for faster 
matrix-matrix multiplication
\cite{Albrecht:2010}.
The fast decomposition of an $n \times m$ matrix with entries in a finite field $\FF_2$ 
is an important issue in algorithmic number theory and cryptanalysis \cite{Shoup:2009,bach:1996}.
In fact, some problems in cryptanalysis and number theory can be transformed in one 
involving a linear system with entries in $\FF_2$. 
The existence of solutions of a linear system can be deduced by analyzing the rank 
of the corresponding matrix.

In this paper we propose a new efficient algorithm to perform 
the matrix factorization for large dense rectangular matrices with entries in $\FF_2$. 
It uses an efficient implementation of M4RM algorithm and uses only row permutations.  
In Section~\ref{sec:3} our notation is given and an appropriate data structure 
to store the matrix is described.
In Section~\ref{sec:4} the non-recursive block decomposition algorithm is presented
and in Section~\ref{sec:5} the corresponding recursive version is described.
Moreover, in Section~\ref{sec:6} some details about the choice of the algorithm parameters are given.
Tests comparing our algorithm with Sage packages~\cite{CGC-misc-sage} are presented in Section~\ref{sec:7}.

\section{The used Matrix Data Structure}\label{sec:3}
Let $\Amat$ be an $n\times m$ matrix having entries in $\FF_2$, i.e. $\Amat \in \matSet(n,m,\FF_2)$,
denoted as 
\begin{EQ}
   \Amat = \pmatrix{
     a_{0,0} & \cdots & a_{0,m-1} \cr
     \vdots & \ddots & \vdots   \cr
     a_{n-1,0} & \cdots & a_{n-1,m-1} \cr
      }.
\end{EQ}
We adopt the following convention for intervals of indices $a..b=\{a,a+1,\ldots,b-1\}$ 
and so, for instance, the sub-matrix $\Amat_{a..b, c..d} \in \matSet(b-a,d-c,\FF_2)$ 
of $\Amat$ represents the intersection of rows of index from $a$ to $b-1$ 
and columns of index from $c$ to $d-1$. The sub-matrix $\Amat_{a..b,\bullet} \in \matSet(b-a,m,\FF_2)$ 
of $\Amat$ is composed by the rows of index from $a$ to $b-1$.
Furthermore, we denote by $\rk(\Amat)$ the rank of $\Amat$, 
i.e. the maximum number of linearly independent row (or column) vectors of $\Amat$. 

We are interested in the computation of the factorization of large dense matrices 
that do not fit into the cache. Thus, a good arrangement of the elements of the matrix in memory 
is important for an efficient data retrieval. Moreover, bits are naturally grouped in words 
whose size is a power of $2$, typically $32$, $64$ or $128$ for larger architectures.
An element of $\FF_2$ is naturally represented as one bit, 
so that elements in $\FF_2$ are naturally grouped in words of $32$, $64$ or $128$ bits.
In particular, a string of elements in $\FF_2$ is packed in an (unsigned) integer.
The advantage of storing multiple elements of $\FF_2$ as an integer is that it guarantees 
a natural parallelism of some operations.

From now on, we denote with $b$ the number of bits of the computer architecture, 
i.e. the number of bits in one machine word. 
Given two integers \verb|x| and \verb|y| whose bits represent elements in $\FF_2$, 
the operation of \emph{exclusive or}, denoted by $\verb|x|\XOR\verb|y|$, is the  
sum $\XOR$ in $\FF_2$ applied to \verb|x| and \verb|y| bitwise;
the \emph{and} operation, denoted by $\verb|x|\AND\verb|y|$, is the  
multiplications $\AND$ in $\FF_2$ applied to \verb|x| and \verb|y| bitwise.
Infact, we have the following formulas
\begin{EQ}
   \verb|x|   = \sum_{i=0}^{b} x_i 2^i,\quad
   \verb|y|   = \sum_{i=0}^{b} y_i 2^i,\quad
   \verb|x|\XOR\verb|y| = \sum_{i=0}^{b} (x_i\XOR y_i) 2^i,\quad
   \verb|x|\AND\verb|y| = \sum_{i=0}^{b} (x_i\AND y_i) 2^i\, .
\end{EQ}
The entries of a matrix are packed into integers that represent groups of 
elements of the matrix itself.
In particular, $b$ consecutive entries in one row are packed into one integer.
The way the integers are arranged changes how one accesses the elements 
of the matrix and you implement the elementary operations.

A matrix $\Amat\in\matSet(n,m,\FF_2)$ is stored into a matrix of non-negative integers
$\AAmat\in\matSet(n,\mu,\NN)$ with $\mu b-b < m \leq \mu b$.
\begin{figure}[!htcb]
\definecolor{ffffff}{rgb}{1,1,1}

\begin{center}
\begin{tikzpicture}[scale=0.71,line cap=round,line join=round,>=triangle 45,x=1.0cm,y=1.0cm]
\clip(-0.49,-1.62) rectangle (5.8,7.68);
\fill[line width=1.2pt,fill=black,fill opacity=0.1] (0,7) -- (0,0.5) -- (4.5,0.5) -- (4.5,7) -- cycle;
\draw (0,0)-- (5,0);
\draw (5,0)-- (5,7);
\draw (5,7)-- (0,7);
\draw (0,7)-- (0,0);
\draw (1,7)-- (1,0);
\draw (2,7)-- (2,0);
\draw (3,7)-- (3,0);
\draw (4,7)-- (4,0);
\draw [line width=1.2pt] (0,7)     -- (0,0.5);
\draw [line width=1.2pt] (0,0.5)   -- (4.5,0.5);
\draw [line width=1.2pt] (4.5,0.5) -- (4.5,7);
\draw [line width=1.2pt] (4.5,7)   -- (0,7);
\draw (0,-0.5) node[anchor=north west] {\parbox{4.52 cm}{$\underbrace{\hspace{3.5cm}}_{\mu b}$}};
\draw (-0.5,0.33) node[anchor=north west] {$\nu$};
\draw (-0.5,0.77) node[anchor=north west] {$n$};
\draw (4.2,7.51) node[anchor=north west] {$m$};
\draw (0,0.1) node[anchor=north west] {$\underbrace{\hspace{0.5cm}}_b$};
\draw (1,0.1) node[anchor=north west] {$\underbrace{\hspace{0.5cm}}_b$};
\draw (2,0.1) node[anchor=north west] {$\underbrace{\hspace{0.5cm}}_b$};
\draw (3,0.1) node[anchor=north west] {$\underbrace{\hspace{0.5cm}}_b$};
\draw (4,0.1) node[anchor=north west] {$\underbrace{\hspace{0.5cm}}_b$};
\end{tikzpicture}
\hbox{
\begin{tikzpicture}[scale=0.71,line cap=round,line join=round,>=triangle 45,x=1.0cm,y=1.0cm]
\clip(-0.49,-1.62) rectangle (5.8,7.68);

\draw (0,0)-- (5,0);
\draw (5,0)-- (5,7);
\draw (5,7)-- (0,7);
\draw (0,7)-- (0,0);
\draw (1,7)-- (1,0);
\draw (2,7)-- (2,0);
\draw (3,7)-- (3,0);
\draw (4,7)-- (4,0);

\fill[line width=1.2pt,fill=black,fill opacity=0.1] (3,5) -- (3,0.5) -- (4.5,0.5) -- (4.5,5) -- cycle;

\fill[pattern color=black,pattern=horizontal lines] (2,5) -- (2,0.5) -- (3,0.5) -- (3,5) -- cycle;

\fill[fill=white] (0,7) -- (2,7)   -- (2,5)   -- (0,5)   -- cycle;
\fill[fill=white] (2,7) -- (4.5,7) -- (4.5,5) -- (2,5)   -- cycle;
\fill[fill=white] (0,5) -- (2,5)   -- (2,0.5) -- (0,0.5) -- cycle;


\draw [line width=1.2pt] (0,7)     -- (0,0.5);
\draw [line width=1.2pt] (0,0.5)   -- (4.5,0.5);
\draw [line width=1.2pt] (4.5,0.5) -- (4.5,7);
\draw [line width=1.2pt] (4.5,7)   -- (0,7);
\draw [line width=1.2pt] (2,7)     -- (2,0.5);
\draw [line width=1.2pt] (2,0.5)   -- (4.5,0.5);
\draw [line width=1.2pt] (4.5,0.5) -- (4.5,5);
\draw [line width=1.2pt] (4.5,5)   -- (0,5);

\draw (0,-0.5) node[anchor=north west] {\parbox{4.52 cm}{$\underbrace{\hspace{3.5cm}}_{\mu b}$}};
\draw (-0.5,0.33) node[anchor=north west] {$\nu$};
\draw (-0.5,0.77) node[anchor=north west] {$n$};
\draw (4.2,7.51) node[anchor=north west] {$m$};
\draw (0,0.1) node[anchor=north west] {$\underbrace{\hspace{0.5cm}}_b$};
\draw (1,0.1) node[anchor=north west] {$\underbrace{\hspace{0.5cm}}_b$};
\draw (2,0.1) node[anchor=north west] {$\underbrace{\hspace{0.5cm}}_b$};
\draw (3,0.1) node[anchor=north west] {$\underbrace{\hspace{0.5cm}}_b$};
\draw (4,0.1) node[anchor=north west] {$\underbrace{\hspace{0.5cm}}_b$};

\draw (1,6)     node[] {\Large$\bm{L}\backslash\bm{U}$};
\draw (1,2.5)   node[] {\Large$\bm{L}$};
\draw (3.25,6)  node[] {\Large$\bm{U}$};
\end{tikzpicture}}
\end{center}

\caption{}
\label{fig:1}
\end{figure}
To access to the element $a_{ij}$ of $\Amat$, we have to determine the corresponding 
column-block $q$ and then the right bit. Precisely, using integer division with remainder
$j = qb+r,\, (0\leq r < b),$
the element $a_{ij}$ corresponds to the $r$-th bit of the integer $\AAscal_{iq}$.
\begin{remark}
  According to the previous notation,
  the least significant bit is on the left respect to the most significant bit,
  i.e. we are using big endian bit order.
\end{remark}
We consider a trivial example.
Let $b=3$ and $\Amat$ be the following $(4 \times 5)$-matrix:
\begin{EQ}
  \Amat = \pmatrix{
     1 & 0 & 1 & 1 & 1  \cr
     1 & 0 & 0 & 0 & 1  \cr
     1 & 1 & 0 & 1 & 0  \cr
     0 & 0 & 1 & 1 & 1  \cr
  }.
\end{EQ}
Since  $m=5$ is not a multiple of $b=3$, we directly memorize the matrix
\begin{EQ}
  \Amat' = \pmatrix{
     1 & 0 & 1 & 1 & 1 & 0 \cr
     1 & 0 & 0 & 0 & 1 & 0 \cr
     1 & 1 & 0 & 1 & 0 & 0 \cr
     0 & 0 & 1 & 1 & 1 & 0 \cr
  }
\end{EQ}
in a $4\times 2$ matrix of $3$ bit integer as follows
\begin{EQ}\label{eq:matrix:A}
  \AAmat = 
  \pmatrix{
     1\cdot2^0+0\cdot2^1+1\cdot2^2 & 1\cdot2^0+1\cdot2^1+0\cdot2^2 \cr
     1\cdot2^0+0\cdot2^1+0\cdot2^2 & 0\cdot2^0+1\cdot2^1+0\cdot2^2 \cr
     1\cdot2^0+1\cdot2^1+0\cdot2^2 & 1\cdot2^0+0\cdot2^1+0\cdot2^2  \cr
     0\cdot2^0+0\cdot2^1+1\cdot2^2 & 1\cdot2^0+1\cdot2^1+0\cdot2^2  \cr
  }
  =
  \pmatrix{
     5 & 3 \cr
     1 & 2 \cr
     3 & 1 \cr
     4 & 3 }.
\end{EQ}
The elements of $\AAmat$ can be organized using 
row-major order or column-major order. 
Our algorithm works better using column-major order. 
For example matrix $\AAmat$ in~\eqref{eq:matrix:A} is stored as
\begin{EQ}
  \pmatrix{ 5 & 1 & 3 & 4 & | & 3 & 2 & 1 & 3}\, .
\end{EQ}
Due to the way the matrix is stored, it is clear that the operations acting on rows 
or block of $b$ consecutive columns (stored as one integer) should be preferred.
Operations acting on isolated columns should be avoided.
In fact, the cost of an operation on a single column on matrix $\Amat$ is about equal (or more) the cost 
of the same operation performed on a group of $b$ columns when the columns are stored 
in a single column of the integer matrix $\AAmat$.

\section{Block decomposition}\label{sec:4}
In this section we give a (non-recursive) block decomposition that holds 
for matrices with entries in any field.  
For simplicity of notation, we are going to describe our strategy for matrices in $\FF_2$.
Let consider $\Amat \in \matSet(n,m,\FF_2)$ such that it can be split in four blocks
\begin{EQ}\label{eq:split:fourblocks}
   \Amat =
   \left[
   \begin{BMAT}(b){c0c}{c0c} \Bmat & \Cmat \\ \Dmat & \Emat \end{BMAT}
   \right],
   \qquad
\end{EQ}
where the block $\Bmat\in\matSet(r,b,\FF_2)$, 
with $r\leq b$, has full rank and the rows of $\Dmat$ are linearly dependent on those of $\Bmat$.
In other words, there exists a matrix $\Ymat$ such that $\Dmat = \Ymat\Bmat$.
The factorization is based on the identity: 
\begin{EQ}\label{eq:factorize:block}
   \left[
   \begin{BMAT}(b){c0c}{c0c} \Bmat & \Cmat \\ \Dmat & \Emat \end{BMAT}
   \right]
   =
   \left[
   \begin{BMAT}(b){c0c}{c0c} \Imat    & \zero \\ \Ymat & \Imat \end{BMAT}
   \right]
   \left[
   \begin{BMAT}(@){c0c}{c0c} \Imat & \zero \\ \zero & \Emat\XOR\Ymat\Cmat \end{BMAT}
   \right]
   \left[
   \begin{BMAT}(@){c0c}{c0c} \Bmat & \Cmat \\ \zero & \Imat \end{BMAT}
   \right].
\end{EQ}
Due to the previous identity, we have that
\begin{EQ}
   \left[
   \begin{BMAT}(@){c0c}{c0c} \Imat & \zero \\ \zero & \Emat\XOR\Ymat\Cmat \end{BMAT}
   \right]
   \left[
   \begin{BMAT}(@){c0c}{c0c} \Bmat & \Cmat \\ \zero & \Imat \end{BMAT}
   \right]
   =
   \left[
   \begin{BMAT}(@){c0c}{c0c} \Bmat & \Cmat \\ \zero & \Emat\XOR\Ymat\Cmat \end{BMAT}
   \right]
\end{EQ}
and notice that $\rk(\Amat)=\rk(\Emat\XOR\Ymat\Cmat)+\rk(\Bmat)=\rk(\Emat\XOR\Ymat\Cmat)+r$.
The sub-matrix $\Emat\XOR\Ymat\Cmat \in \matSet(n-r,m-b,\FF_2)$
is the Schur complement of $\Amat$ (see \cite{haynsworth:1968,zhang:2005}).
Thus, we have reduced the decomposition to a smaller problem.
Applying the same idea as before we can reduce the decomposition  
to smaller and smaller problems with a reduction steps that can be described as follows.\\
Let $\Amat^{(0)}=\Amat$ and $\Pmat^{(0)}$ be a permutation matrix such that 
$\Pmat^{(0)}\Amat^{(0)}$ can be partitioned as
\begin{EQ}
   \Pmat^{(0)}\Amat^{(0)} =
   \left[
   \begin{BMAT}(b){c0c}{c0c} \Bmat^{(0)} & \Cmat^{(0)} \\ \Dmat^{(0)} & \Emat^{(0)} \end{BMAT}
   \right]
   \; \textrm{where}
   \; \Bmat^{(0)}\in\matSet(r_{0},b,\FF_2),
   \quad \rk (\Bmat^{(0)})=r_0\leq b
\end{EQ}
and the rows of $\Dmat^{(0)}$ are linearly dependent 
on those of $\Bmat^{(0)}$, i.e. $\Dmat^{(0)} = \Ymat^{(0)}\Bmat^{(0)}$.
Using the products in~\eqref{eq:factorize:block} we obtain the following
decomposition
\begin{EQ}
  \Pmat^{(0)}\Amat^{(0)}
  =
  \left[
  \begin{BMAT}(b){c0c}{c0c} \Imat & \zero \\ \Ymat^{(0)} & \Imat \end{BMAT}
  \right]
  \left[
  \begin{BMAT}(b){c0c}{c0c} \Imat & \zero \\ \zero & \Amat^{(1)} \end{BMAT}
  \right]
  \left[
  \begin{BMAT}(b){c0c}{c0c} \Bmat^{(0)} & \Cmat^{(0)} \\ \zero & \Imat \end{BMAT}
  \right],
\end{EQ}
where $\Amat^{(1)}=\Emat^{(0)}\XOR\Ymat^{(0)}\Cmat^{(0)} \in \matSet(n-r_{0},m-b,\FF_2)$. 
Now, we can apply the same idea to the Schur complement $\Amat^{(1)}$. 
Let $\Pmat^{(1)}$ be a permutation matrix such that 
\begin{EQ}
   \Pmat^{(1)}
   \Amat^{(1)} =
   \left[
   \begin{BMAT}(b){c0c}{c0c} \Bmat^{(1)} & \Cmat^{(1)} \\ \Dmat^{(1)} & \Emat^{(1)} \end{BMAT}
   \right]
   \; \textrm{where}
   \; \Bmat^{(1)}\in\matSet(r_{1},b,\FF_2),
   \quad 
   \rk (\Bmat^{(1)})=r_1\leq b
\end{EQ}
and $\Dmat^{(1)} = \Ymat^{(1)}\Bmat^{(1)}$.
Due to the decomposition~\eqref{eq:factorize:block}, we reduce further the problem 
of the decomposition of $\Amat^{(1)}$ to the decomposition of 
$\Amat^{(2)}=\Emat^{(1)}\XOR\Ymat^{(1)}\Cmat^{(1)}\in \matSet(n-r_{0}-r_{1},m-2b,\FF_2)$.\\

Going forward, at the $j^{\mathrm{th}}$ stage, the original matrix $\Amat$
is transformed into the sub-matrix $\Amat^{(j)}$; notice that
$\rk(\Amat)=\rk(\Amat^{(j)})+r_{0}+r_{1}+\cdots+r_{j-1}$.
Next, we find a permutation $\Pmat^{(j)}$ such that
$\Pmat^{(j)}\Amat^{(j)}$ can be partitioned as follows
\begin{EQ}\label{eq:permute:Bj}
  \Pmat^{(j)}\Amat^{(j)} =
    \left[
    \begin{BMAT}(b){c0c}{c0c} 
      \Bmat^{(j)} & \Cmat^{(j)} \\ \Dmat^{(j)} & \Emat^{(j)}
    \end{BMAT}
    \right]
    \quad\textrm{with}\quad
    \Bmat^{(j)}\in\matSet(r_{j},b,\FF_2),
\end{EQ}
with $r_{j}\leq b$ and $\Bmat^{(j)}$ full rank and $\Dmat^{(j)} = \Ymat^{(j)}\Bmat^{(j)}$.
Using~\eqref{eq:factorize:block} again, we have 
\begin{EQ}\label{eq:reduction:step}
    \left[
    \begin{BMAT}(@){c0c}{c0c} 
      \Imat & \zero \\ \zero & \Pmat^{(j)}\Amat^{(j)}
    \end{BMAT}
    \right]
   =
   \left[
   \begin{BMAT}(@){c0cc}{c0cc} 
     \Imat & \zero       & \zero \\
     \zero & \Imat       & \zero \\
     \zero & \Ymat^{(j)} & \Imat \end{BMAT}
   \right]
   \left[
     \begin{BMAT}(0){c0cc}{c0cc}
        \Imat & \zero & \zero \\
        \zero & \Imat & \zero \\ 
        \zero & \zero & \Amat^{(j+1)}
      \end{BMAT}
   \right]
   \left[
     \begin{BMAT}(0){c0cc}{c0cc}
        \Imat & \zero & \zero \\
        \zero & \Bmat^{(j)} & \Cmat^{(j)} \\ 
        \zero & \zero & \Imat
      \end{BMAT}
   \right],
\end{EQ}
where
\begin{EQ}\label{eq:Aj}
  \Amat^{(j+1)}=\Emat^{(j)}\oplus \Ymat^{(j)}\Cmat^{(j)}.
\end{EQ}
Note that $\Amat^{(j+1)}$ is the Schur complement of $\Pmat^{(j)}\Amat^{(j)}$ in $\FF_2$
and, after $\mu$ steps, we obtain 
\begin{EQ}\label{eq:pre:final:LU}
  \Pmat\Amat = 
  \underbrace{
  \left[
     \begin{BMAT}(0){c0c}{c0c}
        \Lmat_{11} & \zero \\
        \Lmat_{21} & \Imat
      \end{BMAT}
  \right]
  }_{\Lmat}
  \left[
     \begin{BMAT}(0){c0c}{c0c}
        \Imat & \zero \\
        \zero & \Amat^{(\mu-1)}
      \end{BMAT}
  \right]
  \underbrace{
  \left[
     \begin{BMAT}(0){c0c}{c0c}
        \Umat_{11} & \Umat_{12} \\
        \zero      & \Imat
      \end{BMAT}
  \right]
  }_{\Umat},
\end{EQ}
where $\Lmat$ is non-singular lower triangular matrix and $\Umat$ is the
full rank upper block staircase.
Let $\Pmat^{(\mu-1)}$ be a permutation matrix such that 
$\Pmat^{(\mu-1)}\Amat^{(\mu-1)}$ can be partitioned as
\begin{EQ}
   \Pmat^{(\mu-1)}\Amat^{(\mu-1)} =
   \left[
   \begin{BMAT}(b){c}{c0c} \Bmat^{(\mu-1)} \\ \Dmat^{(\mu-1)}\end{BMAT}
   \right]
   \quad \textrm{where}
   \quad \Bmat^{(\mu-1)}\in\matSet(r_{\mu-1},b',\FF_2),
\end{EQ}
with $b'=m-(\mu-1)b$ satisfying $r_{\mu-1}\leq b'\leq b$
and $\Bmat^{(\mu-1)}$ is full rank. \\
Finally, $\Dmat^{(\mu-1)} = \Ymat^{(\mu-1)}\Bmat^{(\mu-1)}$
and the decomposition~\eqref{eq:pre:final:LU} becomes
\begin{EQ}[rcl]
  \left[
     \begin{BMAT}(0){c0c}{c0c}
        \Imat & \zero \\
        \zero & \Pmat^{(\mu-1)}
      \end{BMAT}
  \right]
  \Pmat\Amat 
  &=&
  \left[
     \begin{BMAT}(0){c0c}{c0c}
        \Lmat_{11} & \zero \\
        \Pmat^{(\mu-1)}\Lmat_{21} & \Imat
      \end{BMAT}
  \right]
  \left[
     \begin{BMAT}(0){c0c}{c0c}
        \Imat & \zero \\
        \zero &
        \left[
        \begin{BMAT}(0){c}{c0c}
          \Imat \\
          \Ymat^{(\mu-1)}
        \end{BMAT}
      \right]
      \Bmat^{(\mu-1)}
      \end{BMAT}
   \right]
  \left[
     \begin{BMAT}(0){c0c}{c0c}
        \Umat_{11} & \Umat_{12} \\
        \zero      & \Imat
      \end{BMAT}
  \right],\\
  &=&
  \underbrace{
  \left[
     \begin{BMAT}(0){c0c}{c0c}
        \Lmat_{11} & \zero \\
        \Pmat^{(\mu-1)}\Lmat_{21} & \left[
        \begin{BMAT}(0){c}{c0c}
          \Imat \\
          \Ymat^{(\mu-1)}
        \end{BMAT}
      \right]
      \end{BMAT}
  \right]
  }_{\Lmat}
  \underbrace{
  \left[
     \begin{BMAT}(0){c0c}{c0c}
        \Umat_{11} & \Umat_{12} \\
        \zero      & \Bmat^{(\mu-1)}
      \end{BMAT}
  \right]
  }_{\Umat},
\end{EQ}
where $\Lmat$ is full rank lower trapezoidal matrix and $\Umat$ is the
full rank upper block staircase.
The previous steps can be resumed in the Lemma:
\begin{lemma}\label{lemma:factorization}
Given any (rectangular) matrix $\Amat$, there exists a permutation $\Pmat$
such that 
\begin{EQ}
   \Pmat\Amat = \Lmat\Umat \, ,
\end{EQ}
where $\Umat$ is full rank upper (block) triangular matrix
and $\Lmat$ is a full rank lower trapezoidal matrix.
\end{lemma}

Clearly, $\rk(\Amat)$ is the rank of the matrix $\Umat$ which is
the number of its rows: $\rk(\Amat)=\sum_{j=0}^{\mu-1} r_{j}$.

Observe that the computation at the step $j$ involves
matrix-matrix multiplications to obtain the Schur complement~\eqref{eq:Aj}
which are the most costly operations of the presented algorithm.
The multiplication of a $p \times b$-matrix by a $b \times q$-matrix, in case $p\gg b$, 
can be efficiently performed using the M4RM algorithm 
\cite{Arlazarov:1970,aho:1974,Albrecht:2010}. Thus, 
in the computation of $\Amat^{(j+1)}$ it is convenient to use the M4RM algorithm,
because this block operation is more efficient 
than the usual row operations.

This decomposition is based on the selection of $\Bmat^{(j)}$ and the construction of 
$\Ymat^{(j)}$ at each step. An efficient algorithm for this will be discussed
in the next sections.
For this purpose the incremental construction of an inverse of
a matrix and pseudo-inverse construction is a necessary tool.

\begin{remark}
The decomposition described in Lemma~\ref{lemma:factorization}
when $b=1$ is equivalent to the PLE factorization described
in~\cite{Albrecht:2011,Jeannerod:2011}.
However, computation with $b=1$ is not convenient losing natural 
parallelism of integer operations.
\end{remark}

\subsection{The computation of $\Ymat^{(j)}$}\label{Comput:W}
Assume that the permutation $\Pmat$ applied to matrix $\Amat$ results in
\begin{EQ}\label{eq:permute:Bj:low}
   \Pmat\Amat
   =
   \left[
    \begin{BMAT}(b){c0c}{c0c} 
      \Bmat & \Cmat \\ \Dmat & \Emat
    \end{BMAT}
   \right],
\end{EQ}
where the rectangular matrix $\Bmat\in\matSet(r,b,\FF_2)$, with $r\leq b$, has full rank 
and satisfies $\Dmat=\Ymat\Bmat$ for an opportune matrix $\Ymat$ which we have to determine.

In case $r=b$ matrix $\Bmat$ is non-singular and we easily deduce that $\Ymat=\Dmat \Bmat^{-1}$.
Instead, when the matrix $\Bmat$ is full rank with less rows than columns, 
a \emph{pseudo-inverse} $\Bmat^\dagger$ has to be computed.
For example the pseudo-inverse of Moore-Penrose~\cite{Israel:2004} is given by 
\begin{EQ}
   \Bmat^\dagger = \Bmat^T(\Bmat\Bmat^T)^{-1}
\end{EQ}
and satisfies
\begin{EQ}[rcl]
   \Bmat\Bmat^\dagger &=& \Bmat\Bmat^T(\Bmat\Bmat^T)^{-1} = \Imat, \\
   \Dmat\Bmat^\dagger &=& \Ymat\Bmat\Bmat^\dagger = \Ymat.
\end{EQ}
Thus, $\Ymat$ can be computed by simple right multiplication by
a pseudo-inverse of $\Bmat$.\\
Although the use of the Moore-Penrose's pseudo-inverse
is correct, a more efficient pseudo-inverse can be constructed.
Let $\Jmat\in\matSet(b,r,\FF_2)$ be an \emph{insertion} matrix 
whose effect is to insert $b-r$ zero-rows into a matrix with $r$ rows.
Let $\Rmat\in\matSet(b,b,\FF_2)$ be the matrix containing
linearly independent rows which makes the square matrix $\Jmat\Bmat+\Rmat$ non-singular.
Notice that $\Jmat^T$ is a \emph{projection} matrix which satisfies
\begin{EQ}\label{eq:JR:prop}
  \Jmat^T\Jmat=\Imat,\qquad
  \Jmat^T\Rmat=\zero.
\end{EQ}
If $r=b$, i.e. if $\Bmat$ is square and non-singular, insertions are not necessary and we have $\Jmat=\Imat$, $\Rmat=\zero$.
\begin{example}
  In case $b=3$ and $r=2$ we can see a situation as
  \begin{EQ}
    \Jmat_j = 
    \left[
      \begin{BMAT}(b){cc}{ccc} 1 & 0 \\ 0 & 0 \\ 0 & 1 \end{BMAT}
    \right],
    \qquad
    \Rmat_j = 
    \left[
      \begin{BMAT}(b){ccc}{ccc} 0 & 0 & 0 \\ 0 & 1 & 0 \\ 0 & 0 & 0 \end{BMAT}
    \right].
  \end{EQ}
\end{example}
Let consider $\Zmat=(\Jmat\Bmat+\Rmat)^{-1}\Jmat$. Due to the properties in~\eqref{eq:JR:prop} and the relation $\Dmat=\Ymat\Bmat$, we have
\begin{EQ}[rcl]\label{eq:pre:W}
   \Bmat\Zmat
   &=&\Bmat(\Jmat\Bmat+\Rmat)^{-1}\Jmat
   =\Jmat^T\Jmat\Bmat(\Jmat\Bmat+\Rmat)^{-1}\Jmat\\
   &=&\Jmat^T(\Jmat\Bmat+\Rmat)(\Jmat\Bmat+\Rmat)^{-1}\Jmat
   =\Jmat^T\Jmat=\Imat\\
   \Dmat\Zmat &=&\Ymat\Bmat\Zmat =\Ymat \, .
\end{EQ}
Thus the matrix $\Zmat$ has the same effect of the pseudo-inverse (it is a pseudo-inverse
different from the Moore-Penrose one) and it is used in the computation
of $\Ymat$ during the factorization procedure.
In order to build $\Zmat$, we need an algorithm to compute the inverse of 
a (small) square matrix in $\FF_2$; it will be treated in the next section.

\subsection{Incremental construction of the inverse of a square matrix in $\FF_2$}
\label{sec:incremental}

Let $\Bmat$ non-singular with all principal minors
non-singular, it is possible to incrementally build its inverse $\Zmat$.
This requirement is not restrictive because every non-singular matrix by a row permutation satisfies it 
(due to Gauss $LU$ decomposition, see~\cite{kincaid:2002} page $156$ Theorem $1$). 
Let $\Bmat_k$ be the $k^\mathrm{th}$ principal minor of $\Bmat$
and $\Zmat_k$ its inverse. We can directly obtain $\Bmat_{k+1}^{-1}$ from $\Bmat_{k}^{-1}$
using the following factorization
\begin{EQ}\label{eq:Bk:FACT}
   \Bmat_{k+1} =
   \left[
     \begin{BMAT}(b){cc}{cc}
       \Bmat_k & \cvec \\ 
       \rvec^T & \alpha
     \end{BMAT}
   \right]
   =
   \left[
     \begin{BMAT}(b){cc}{cc}
       \Imat & \zero \\ 
       \rvec^T\Zmat_k & 1
     \end{BMAT}
   \right]
   \left[
     \begin{BMAT}(b){cc}{cc}
       \Bmat_k & \cvec \\ 
             0 & 1
     \end{BMAT}
   \right]
\end{EQ}
which holds when
\begin{EQ}\label{eq:cond:nonsing:Bk}
  \alpha \oplus \rvec^T\Zmat_k\cvec=1.
\end{EQ}
\begin{remark}
In case of $\FF_q$, the previous condition becomes $\alpha - \rvec^T\Zmat_k\cvec=\beta$ where 
$\alpha, \beta \in \FF_q$ and $\beta\not=0$. Moreover, the last block in~\eqref{eq:Bk:FACT} becomes
\begin{EQ}
   \left[
     \begin{BMAT}(b){cc}{cc}
       \Bmat_k & \cvec \\ 
             0 & \beta
     \end{BMAT}
   \right].
\end{EQ}
\end{remark}
Inverting factorization~\eqref{eq:Bk:FACT} we get immediately:
\begin{EQ}\label{eq:Bk_k1}
   \Bmat_{k+1}^{-1} =
   \left[
     \begin{BMAT}(b){cc}{cc}
       \Bmat_k & \cvec \\ 
       \rvec^T & \alpha
     \end{BMAT}
   \right]^{-1}
   =
   \left[
     \begin{BMAT}(r){cc}{cc}
       \Zmat_k & \Zmat_k\cvec \\ 
       \zero^T & 1
     \end{BMAT}
   \right]
   \left[
     \begin{BMAT}(r){cc}{cc}
       \Imat & \zero \\ 
       \rvec^T\Zmat_k & 1
     \end{BMAT}
   \right]
\end{EQ}
(Notice that $\Bmat_1$ is a $1\times 1$ matrix and thus $\Bmat_1=\Zmat_1=[1]$).
Therefore, due to~\eqref{eq:Bk_k1} we can compute $\Zmat_{k+1}$ from $\Zmat_{k}$ as  
\begin{EQ}\label{eq:fact:pre:Z}
   \Zmat_{k+1} =\Bmat_{k+1}^{-1} 
   =
    \left[
     \begin{BMAT}(r){cc}{cc}
       \Zmat_k\XOR(\Zmat_k\cvec)(\rvec^T\Zmat_k) & \Zmat_k\cvec \\
       \rvec^T\Zmat_k & 1
     \end{BMAT}
   \right],
\end{EQ}
which needs two matrix-vector multiplication and a rank one update. 
Moreover, setting $\widetilde{\cvec}=\Zmat_k\cvec$, the last matrix in~\eqref{eq:fact:pre:Z}
can be written as a matrix-matrix product:
\begin{EQ}\label{eq:fact:pre:Z1}
    \Zmat_{k+1} =
    \left[
     \begin{BMAT}(r){cc}{cc}
       \Zmat_k\XOR\widetilde{\cvec}(\rvec^T\Zmat_k) & \widetilde{\cvec}\\
       \rvec^T\Zmat_k & 1
     \end{BMAT}
   \right]
   = 
  \Hmat_{k+1}
    \left[
     \begin{BMAT}(r){cc}{cc}
       \Zmat_k & \zero \\
       \zero & 1
     \end{BMAT}
   \right], \quad
     \Hmat_{k+1}=
    \left[
     \begin{BMAT}(r){cc}{cc}
       \Imat\XOR\widetilde{\cvec}\rvec^T & \widetilde{\cvec} \\
       \rvec^T & 1 \\
     \end{BMAT}
   \right].
   \qquad
\end{EQ}
Let $\Mmat_k$ be the matrix obtained multiplying by $\Zmat_k$ the first $k$ rows of $\Bmat$:
\begin{EQ}\label{eq:Mk}
   \Mmat_k =
   \left[
     \begin{BMAT}(r){ccc}{ccc}
       \Imat   & \Zmat_k\cvec & \Zmat_k\Cmat \\ 
       \rvec^T & \alpha       & \evec^T \\
       \Dmat   & \dvec       & \Emat
     \end{BMAT}
   \right]
   =
   \left[
     \begin{BMAT}(r){ccc}{ccc}
       \Imat   & \widetilde{\cvec} & \widetilde{\Cmat} \\ 
       \rvec^T & \alpha       & \evec^T \\
       \Dmat   & \dvec       & \Emat
     \end{BMAT}
   \right]
   \quad
   \textrm{where}
   \quad
   \Bmat
   = 
   \left[
     \begin{BMAT}(r){ccc}{ccc}
       \Bmat_k  & \cvec  & \Cmat   \\ 
       \rvec^T  & \alpha & \evec^T \\
       \Dmat    & \dvec  & \Emat
     \end{BMAT}
   \right].
\end{EQ}
Due to~\eqref{eq:fact:pre:Z1}, the update of $\Mmat_{k+1}$ results in
\begin{EQ}[rcl]
    \Mmat_{k+1} &=&
       \left[
     \begin{BMAT}(r){cc}{cc}
       \Zmat_{k+1} & \zero \\ 
       \zero        & \Imat \\
     \end{BMAT}
   \right]
   \Bmat=
   \left[
     \begin{BMAT}(r){c.c}{c.c}
       \Hmat_{k+1} & \zero \\
       \zero & \Imat
     \end{BMAT}
   \right]
    \left[
     \begin{BMAT}(r){cc.c}{cc.c}
       \Zmat_k & \zero & \zero \\
       \zero & 1 & \zero \\
       \zero & \zero & \Imat
     \end{BMAT}
   \right]\Bmat
   =
    \left[
     \begin{BMAT}(r){c.c}{c.c}
       \Hmat_{k+1} & \zero \\
       \zero & \Imat
     \end{BMAT}
   \right]
   \Mmat_k 
\end{EQ}
and thus, $\Hmat_{k+1}$ is used in updating both $\Mmat_{k+1}$ and $\Zmat_{k+1}$.
Moreover, we obtain the following update formulas:
\begin{EQ}
   \Hmat_{k+1}
   \left[
     \begin{BMAT}(r){c}{cc}
       \widetilde{\Cmat}   \\
       \evec^T \\
     \end{BMAT}
   \right]
   =
   \left[
     \begin{BMAT}(r){cc}{cc}
       \Imat\XOR\widetilde{\cvec}\rvec^T & \widetilde{\cvec} \\
       \rvec^T & 1 \\
     \end{BMAT}
   \right]
   \left[
     \begin{BMAT}(r){c}{cc}
       \widetilde{\Cmat}   \\
       \evec^T \\
     \end{BMAT}
   \right]
   =
   \left[
     \begin{BMAT}(r){c}{cc}
       \widetilde{\Cmat}\XOR\widetilde{\cvec}\svec^T \\
       \svec^T \\
     \end{BMAT}
   \right]
   \quad
   \textrm{where}
   \quad
   \svec^T = \rvec^T\widetilde{\Cmat}\XOR\evec^T .
\end{EQ}
Note that it is possible to perform the incremental build when all principal minors are non-singular, 
which is equivalent to satisfy equation~\eqref{eq:cond:nonsing:Bk} for all $k$.
Such a condition is used to dynamically select linear independent rows in the construction of $\Bmat^{(j)}$.

\subsection{Construction of row permutation}\label{sec:permutation}
The steps of reduction from $\Amat^{(j)}$ to $\Amat^{(j+1)}$ need the computation 
of permutation matrix $\Pmat^{(j)}$ that reorders the rows of $\Amat^{(j)}$ 
in order to obtain $\Pmat^{(j)}\Amat^{(j)}$ partitioned as in \eqref{eq:permute:Bj}.
The block $\Bmat^{(j)}$ must be non-singular with all principal minors non-singular.
This permutation can be computed using only the first $b$ columns of $\Amat^{(j)}$ 
and must satisfy 
\begin{EQ}
  \Pmat^{(j)}\Amat^{(j)}_{\bullet, 0..b}=
  \left[
    \begin{BMAT}(b){c}{c0c} 
      \Bmat^{(j)} \\ \Dmat^{(j)}
    \end{BMAT}
  \right].
\end{EQ}
The selection of the rows and the construction of the inverse of $\Bmat^{(j)}$ 
are done together using incremental update of Section~\ref{sec:incremental},
where the selected rows must satisfy condition~\eqref{eq:cond:nonsing:Bk}.
Observe that to check condition~\eqref{eq:cond:nonsing:Bk}, the $\ell$-row of the sub-matrix $\Amat^{(j)}_{\bullet, 0..b}$ is partitioned as
$\Amat^{(j)}_{\ell, 0..b}=(\rvec^T,\alpha,\bm{\beta}^T)$ 
where the $\bm{\beta}^T$ portion of the row is ignored in the computation.

\subsection{Insertion of linearly independent rows}\label{sec:insertion}
In the computation of row permutation it may happen that condition~\eqref{eq:cond:nonsing:Bk} 
is not satisfied by all the rows of the column block $\Amat^{(j)}_{\bullet, 0..b}$.
Obviously, such a situation does not arise if $n\geq m$ and when the matrix is full rank.\\ 
At this stage, standard algorithms introduce column permutations
to satisfy condition~\eqref{eq:cond:nonsing:Bk}.
If such a condition cannot be satisfied even using column permutations, 
it means that last rows are linearly dependent on the previous ones
and then the algorithm must end.
Column permutations are not executed in our algorithm,
so that new linearly independent rows are inserted using the two matrices $\Jmat$ and $\Rmat$,
 introduced in Section~\ref{Comput:W}.
We observe that it is easy to build a row that satisfies condition~\eqref{eq:cond:nonsing:Bk} and that, 
in particular, the row $(\rvec^T,\alpha,\bm{\beta}^T)=(\zero^T,1,\zero^T)$
trivially satisfies it.
In practice, the presented process mixes row permutations and row insertions; 
it can be respectively split in a row permutations followed by a rows insertions:
\begin{EQ}
  \Pmat^{(j)}\Amat^{(j)}_{\bullet, 0..b}=
  \left[
    \begin{BMAT}(b){c}{c0c} 
      \Bmat^{(j)} \\ \Dmat^{(j)}
    \end{BMAT}
  \right],
  \qquad
  \Jmat_j\Bmat^{(j)}+\Rmat_j.
\end{EQ}
Notice that the permutation $\Pmat^{(j)}$ and the insertion $\Jmat_j$ are chosen in such a way the 
square block $\Jmat_j\Bmat^{(j)}+\Rmat_j$ has all principal minors non-singular and satisfies~\eqref{eq:JR:prop}.\\

The following procedure performs the operation described in 
Sections~\ref{sec:incremental}-\ref{sec:permutation} and \ref{sec:insertion}
to obtain an incremental construction of the pseudo-inverse $\Zmat$. 

\begin{small}
\begin{procedure}[H]
  \caption{buildZ($\AAscal,r$)}
  \label{alg:buildZ}
  \SetKwFunction{popCount}{popCount}
  \SetKwFunction{searchRow}{searchRow}
  \SetKwFunction{swapRows}{swapRows}
  \SetKwFunction{columnExtract}{columnExtract}
  \SetKw{goto}{go to line}
  \SetKw{Break}{break}

  $i_b\assign 0$\;
  \lFor{$j=0..b$}{
    $\ZZscal_j\assign 2^j$;\quad
    $\MMscal_j\assign 2^j$;\quad
    $\PPscal_j\assign -1$\;
  }
  \For{$k_{\mathrm{bit}}=0..b$}{%
    $c\assign 2^{k_{\mathrm{bit}}}$\;
    \For{$i=0..k_{\mathrm{bit}}$}{%
      \lIf{$\MMscal_{i}\AND 2^{k_{\mathrm{bit}}} \neq 0$}{$c\assign c \XOR 2^i$\;}
    }
    \For{$i=i_b..r$}{
      \If{\popCount{$c\AND\AAscal_{i}$} \emph{is odd}}{
        $\PPscal_{k_{\mathrm{bit}}}\assign i$;\quad
        $\AAscal_{i}\rightleftharpoons\AAscal_{i_b}$;\quad
        $\MMscal_{k_{\mathrm{bit}}}\assign\AAscal_{i_b}$;\quad
        $i_b\assign i_b+1$;\quad
        \Break ;
      }
    }
    $y\assign \MMscal_{k_{\mathrm{bit}}}$;\quad
    \tcp{Update of $\ZZscal$ and $\MMscal$}
    \For{$j=0..k_{\mathrm{bit}}$}{
      \lIf{$y\AND 2^j\neq 0$}{
        $\ZZscal_{k_{\mathrm{bit}}}\assign \ZZscal_{k_{\mathrm{bit}}}\XOR\ZZscal_{j}$;\quad
        $\MMscal_{k_{\mathrm{bit}}}\assign \MMscal_{k_{\mathrm{bit}}}\XOR\MMscal_{j}$\;
      }
    }
    \For{$j=0..k_{\mathrm{bit}}$}{
      \lIf{$c\AND 2^j\neq 0$}{
        $\ZZscal_{j}\assign \ZZscal_{k_{\mathrm{bit}}}\XOR\ZZscal_{j}$;\quad
        $\MMscal_{j}\assign \MMscal_{k_{\mathrm{bit}}}\XOR\MMscal_{j}$\;
      }
    }
  }
  $m\assign 0$\;
  \For{$i=0..b$}{
    \uIf{$\PPscal_i=-1$}{
      \tcp{$\overline{m}$ is the integer with bits complemented}
      \lFor{$j=0..b$}{$\ZZscal_{j}\assign (\ZZscal_{j}/2\AND \overline{m})\XOR(\ZZscal_{j}\AND m)$\;}
    }
    \Else{
      $m\assign 2m+1$\;
    }
  }
\end{procedure}
\end{small}

\section{Block recursive algorithm}\label{sec:5}
In this section we describe a recursive version of the algorithm presented in Section~\ref{sec:4}. 
Let consider $\Amat \in \matSet(n,m,\FF_2)$ split into two sub-matrices $\Amat_L\in \matSet(n,p,\FF_2)$ and 
$\Amat_R\in \matSet(n,m-p,\FF_2)$ with $p=[m/2]$ such that
\begin{EQ}
   \Amat
   =
   \left[
    \begin{BMAT}(b){c0c}{c} 
      \Amat_L & \Amat_R
    \end{BMAT}
   \right].
\end{EQ}
Applying the decomposition in Lemma 4.1 to the left part $\Amat_L$, we have $\Pmat\Amat_L = \Lmat\Umat $
and so
\begin{EQ}
   \Pmat \Amat
   =
   \left[
    \begin{BMAT}(@){c0c}{c} 
      \Pmat\Amat_L & \Pmat \Amat_R 
    \end{BMAT}
   \right]
   = 
   \left[
    \begin{BMAT}(@){c0c}{c} 
      \Lmat
      \Umat & 
      \Pmat \Amat_R
    \end{BMAT}
   \right] \, .
\end{EQ}
Since 
$\Pmat\Amat_R  = \left[\begin{BMAT}(@){c}{c0c} \Cmat \\ \Dmat \end{BMAT}\right]$, 
where  $\Cmat \in \matSet(r,m-p,\FF_2)$ and $\Dmat \in \matSet(n-r,m-p,\FF_2)$, we obtain
\begin{EQ}[rcl]
  \Pmat \Amat = 
   \left[
    \begin{BMAT}(@){c0c}{c} 
      \left(
   \begin{BMAT}(@){c}{cc} \Lmat_{0..r,\bullet} \\ \Lmat_{r..n,\bullet}\end{BMAT}
   \right)
   \Umat & 
    \left(
   \begin{BMAT}(@){c}{cc} \Cmat \\ \Dmat\end{BMAT}
   \right)    \end{BMAT}
   \right]
  =
   \left[
    \begin{BMAT}(b){c0c}{c0c} 
      \Lmat_{0..r,\bullet}  & \zero \\ \Lmat_{r..n,\bullet} &\Imat  
    \end{BMAT}
   \right]
   \left[
    \begin{BMAT}(@){c0c}{c0c} 
    \Umat & \Lmat_{0..r,\bullet}^{-1} \Cmat \\
    \zero & \Dmat\XOR\Lmat_{r..n,\bullet} \Lmat_{0..r,\bullet}^{-1} \Cmat  
    \end{BMAT}
   \right]\, .
\end{EQ}
Let $\Amat' = \Dmat\XOR\Lmat_{r..n,\bullet} \Lmat_{0..r,\bullet}^{-1} \Cmat$, we can recursively apply the factorization
$\Pmat'\Amat' = \Lmat'\Umat'$ and we obtain
\begin{EQ}[rcl]
  \Pmat \Amat 
  = 
   \left[
    \begin{BMAT}(b){c0c}{c0c} 
      \Lmat_{0..r,\bullet}  & \zero \\ \Lmat_{r..n,\bullet} &\Imat  
    \end{BMAT}
   \right]
   \left[
    \begin{BMAT}(@){c0c}{c0c} 
    \Umat & \Lmat_{0..r,\bullet}^{-1}\Cmat \\
    \zero & \Amat' 
    \end{BMAT}
   \right]
  = 
   \left[
    \begin{BMAT}(b){c0c}{c0c} 
      \Lmat_{0..r,\bullet}  & \zero \\ \Lmat_{r..n,\bullet} & (\Pmat')^{-1}  
    \end{BMAT}
   \right]
   \left[
    \begin{BMAT}(@){c0c}{c0c} 
    \Umat & \Lmat_{0..r,\bullet}^{-1}\Cmat \\
    \zero & \Lmat'\Umat' 
    \end{BMAT}
   \right]
\end{EQ}
and then
\begin{EQ}[rcl]
   \left[
    \begin{BMAT}(b){c0c}{c0c} 
      \Imat  & \zero \\ \zero & \Pmat'
    \end{BMAT}
   \right]
  \Pmat \Amat 
  &=& 
   \left[
    \begin{BMAT}(b){c0c}{c0c} 
      \Lmat_{0..r,\bullet}  & \zero \\ \Pmat'\Lmat_{r..n,\bullet} & \Lmat'  
    \end{BMAT}
   \right]
   \left[
    \begin{BMAT}(@){c0c}{c0c} 
    \Umat & \Lmat_{0..r,\bullet}^{-1}\Cmat \\
    \zero & \Umat' 
    \end{BMAT}
   \right].\\
\end{EQ}
\begin{remark}
Notice that here we perform the matrix-matrix multiplications using our implementation of the Strassen's algorithm.
\end{remark}

\section{Performance tuning}\label{sec:6}
The algorithm presented in Section~\ref{sec:4} completely avoids column permutations 
and the cost of the decomposition is given by the cost of the matrix-matrix multiplication.
The cost to multiply $\Amat \in \matSet(n,b,\FF_2)$ and $\Bmat \in \matSet(b,b,\FF_2)$,
using the M4RM algorithm, is approximately given by the costs of the $\Xor$ operations.
In particular, if we neglect the cost of the memory access and other minor costs, we obtain that
\begin{EQ}\label{eq:M4RM:cost:pre}
  \textrm{cost of M4RM}=
  \texttt{T}_\ell^{b} + n \texttt{R}_\ell^{b},
  \qquad
  \texttt{T}_\ell^{b}=\sum_{j=1}^{K}(2^{\ell_j}-1)\Xor_{b},
  \quad
  \texttt{R}_\ell^{b} = K\,\Xor_{b}\qquad
\end{EQ}
where
\begin{enumerate}
\item $\ell_j$ are positive integers such that $\ell_1+\ell_2+\cdots+\ell_{K}=b$;
\item $2^{\ell_j}$ is the size of the $j^{th}$ of the $K$ tables used to perform M4RM algorithm; 
\item $\texttt{T}_\ell^{b}$ is the cost of the tables construction;
\item $\texttt{R}_\ell^{b}$ is the cost of the rows operations;
\item $\Xor_{b}$ is the cost of the $\Xor$ operation using integer of $b$ bit size.
\end{enumerate}
When $K$ divides $b$, we assume $\ell_j=c=b/K$ and the expression in~\eqref{eq:M4RM:cost:pre}
becomes
\begin{EQ}\label{eq:M4RM:cost}
  \textrm{cost of M4RM}=
  \texttt{T}_c^{b} + n \texttt{R}_c^{b},
  \qquad
  \texttt{T}_c^{b}=\left(\frac{b}{c}\right)(2^{c}-1)\Xor_{b},
  \quad
  \texttt{R}_c^{b} = \left(\frac{b}{c}\right)\Xor_{b}\, .
\end{EQ}
In order to simplify the cost estimation, from now on, we assume that $K$ divides $b$
and we perform the cost analysis starting by~\eqref{eq:M4RM:cost}.
In Subsection~\ref{sec:size of b and c} we are going to discuss about the choices of the parameters $b$ and $c$. 
Moreover, in Subsection~\ref{sec:switching point} we will give an idea about the switching point from the M4RM multiplication and the Strassen multiplication. Finally, in Subsection~\ref{sec:complexity estimation} a complexity estimation of our algorithm is given.

\subsection{Integer word-size and performance of M4RM}
\label{sec:size of b and c}
Let us consider two matrices $\Amat \in \matSet(n,2b,\FF_2)$ and $\Bmat\in\matSet(2b,2b,\FF_2)$. 
Due to~\eqref{eq:M4RM:cost} and choosing $2b$ bits as word-size, the cost to perform the M4RM algorithm to multiply $\Amat$ and $\Bmat$ is approximately $\texttt{T}_c^{2b} + n \texttt{R}_c^{2b}$. On the other hand, using a word-size of $b$ bits, the cost of M4RM becomes
$4\texttt{T}_c^{b} + n 4 \texttt{R}_c^{b}$.\\
Considering the cost of the rows contribution and the respective ratio, we have that
\begin{EQ}
  \frac{\textrm{Cost of $b$-bits row}}{\textrm{Cost of $2b$-bits row}}
  =\frac{4n\texttt{R}_c^{b}}{n\texttt{R}_{c}^{2b}}
  =2\left(\frac{\Xor_{b}}{\Xor_{2b}}\right)
\end{EQ}
and so, in case of $\Xor_{2b}<2\Xor_{b}$, it is convenient to use $2b$ bits.\\
The cost of the tables construction is respectively given by $4\texttt{T}_c^{b}$ and $\texttt{T}_{c}^{2b}$
and the corresponding ratio is then
\begin{EQ}
  \frac{\textrm{Cost per $b$-bits tables}}{\textrm{Cost per $2b$-bits tables}}
  =\frac{4\texttt{T}_c^{b}}{\texttt{T}_{c}^{2b}}
  =2\left(\frac{{\Xor}_{b}}{\Xor_{2b}}\right).
\end{EQ}
Since the cost of the $\Xor_{32} \approx \Xor_{64}$ in $64$ bit architecture,
it is convenient to choose $b=64$. For hardware that supports SSE $128$ bit instructions,
the following relation $2\Xor_{64}>\Xor_{128}\geq\Xor_{64}$ holds 
and thus it is convenient to use $b=128$.
\begin{table}
  \caption{Cost of tables construction and cost of the rows operation for M4RM.}
  \label{tab1:rows and tables costs}
  \begin{center}
    \begin{tikzpicture}[scale=0.4, line cap=round,line join=round,>=triangle 45,x=1.0cm,y=1.0cm]
      \clip(-4.3,-1.54) rectangle (3.4,6.8);
      \fill[fill=black,pattern=horizontal lines] (-4,4.84) -- (-4,-1.16) -- (-2.5,-1.16) -- (-2.5,4.84) -- cycle;
      \fill[fill=black,pattern=horizontal lines] (1.5,4.84) -- (1.5,3.34) -- (3,3.34) -- (3,4.84) -- cycle;
      \fill[fill=black,pattern=horizontal lines] (-0.5,4.84) -- (-0.5,-1.16) -- (1,-1.16) -- (1,4.84) -- cycle;
      \draw [line width=1pt] (-4,4.84)-- (-4,-1.16);
      \draw [line width=1pt] (-4,-1.16)-- (-2.5,-1.16);
      \draw [line width=1pt] (-2.5,-1.16)-- (-2.5,4.84);
      \draw [line width=1pt] (-2.5,4.84)-- (-4,4.84);
      \draw [line width=1pt] (1.5,4.84)-- (1.5,3.34);
      \draw [line width=1pt] (1.5,3.34)-- (3,3.34);
      \draw [line width=1pt] (3,3.34)-- (3,4.84);
      \draw [line width=1pt] (3,4.84)-- (1.5,4.84);
      \draw [line width=1pt] (-0.5,4.84)-- (-0.5,-1.16);
      \draw [line width=1pt] (-0.5,-1.16)-- (1,-1.16);
      \draw [line width=1pt] (1,-1.16)-- (1,4.84);
      \draw [line width=1pt] (1,4.84)-- (-0.5,4.84);
      \draw [line width=1pt] (-0.4,6.32) node[anchor=north west] {\Large$\bm{A}$};
      \draw [line width=1pt] ( 1.6,6.32) node[anchor=north west] {\Large$\bm{B}$};
      \draw [line width=1pt] (-4.0,6.32) node[anchor=north west] {\Large$\bm{C}$};
      \draw [line width=1pt] (-2.24,3.18) node[anchor=north west] {\Large$\bm{=}$};
    \end{tikzpicture}
    \qquad\qquad\qquad
    \begin{tikzpicture}[scale=0.4, line cap=round,line join=round,>=triangle 45,x=1.0cm,y=1.0cm]
      \clip(-4.74,-1.74) rectangle (6.72,5.86);
      \fill[fill=black,pattern=horizontal lines] (-3,4.5) -- (-3,-1.5) -- (-1.5,-1.5) -- (-1.5,4.5) -- cycle;
      \fill[fill=black,pattern=horizontal lines] (3.5,4.5) -- (3.5,3) -- (5,3) -- (5,4.5) -- cycle;
      \fill[fill=black,pattern=horizontal lines] (0,4.5) -- (0,-1.5) -- (1.5,-1.5) -- (1.5,4.5) -- cycle;
      \fill[fill=black,pattern=horizontal lines] (5,4.5) -- (5,3) -- (6.5,3) -- (6.5,4.5) -- cycle;
      \fill[fill=black,pattern=horizontal lines] (3.5,3) -- (3.5,1.5) -- (5,1.5) -- (5,3) -- cycle;
      \fill[fill=black,pattern=horizontal lines] (5,3) -- (5,1.5) -- (6.5,1.5) -- (6.5,3) -- cycle;
      \fill[fill=black,pattern=horizontal lines] (1.5,4.5) -- (1.5,-1.5) -- (3,-1.5) -- (3,4.5) -- cycle;
      \fill[fill=black,pattern=horizontal lines] (-4.5,4.5) -- (-4.5,-1.5) -- (-3,-1.5) -- (-3,4.5) -- cycle;
      \draw [line width=1pt] (-3,4.5)-- (-3,-1.5);
      \draw [line width=1pt] (-3,-1.5)-- (-1.5,-1.5);
      \draw [line width=1pt] (-1.5,-1.5)-- (-1.5,4.5);
      \draw [line width=1pt] (-1.5,4.5)-- (-3,4.5);
      \draw [line width=1pt] (3.5,4.5)-- (3.5,3);
      \draw [line width=1pt] (3.5,3)-- (5,3);
      \draw [line width=1pt] (5,3)-- (5,4.5);
      \draw [line width=1pt] (5,4.5)-- (3.5,4.5);
      \draw [line width=1pt] (0,4.5)-- (0,-1.5);
      \draw [line width=1pt] (0,-1.5)-- (1.5,-1.5);
      \draw [line width=1pt] (1.5,-1.5)-- (1.5,4.5);
      \draw [line width=1pt] (1.5,4.5)-- (0,4.5);
      \draw (0.8,6.) node[anchor=north west] {\Large$\bm{A}$};
      \draw (4.3,6.) node[anchor=north west] {\Large$\bm{B}$};
      \draw (-3.7,6.) node[anchor=north west] {\Large$\bm{C}$};
      \draw (-1.5,3) node[anchor=north west] {\Large$\bm{=}$};
      \draw [line width=1pt] (5,4.5)-- (5,3);
      \draw [line width=1pt] (5,3)-- (6.5,3);
      \draw [line width=1pt] (6.5,3)-- (6.5,4.5);
      \draw [line width=1pt] (6.5,4.5)-- (5,4.5);
      \draw [line width=1pt] (3.5,3)-- (3.5,1.5);
      \draw [line width=1pt] (3.5,1.5)-- (5,1.5);
      \draw [line width=1pt] (5,1.5)-- (5,3);
      \draw (5,3)-- (3.5,3);
      \draw [line width=1pt] (5,3)-- (5,1.5);
      \draw [line width=1pt] (5,1.5)-- (6.5,1.5);
      \draw [line width=1pt] (6.5,1.5)-- (6.5,3);
      \draw (6.5,3)-- (5,3);
      \draw (1.5,4.5)-- (1.5,-1.5);
      \draw [line width=1pt] (1.5,-1.5)-- (3,-1.5);
      \draw [line width=1pt] (3,-1.5)-- (3,4.5);
      \draw [line width=1pt] (3,4.5)-- (1.5,4.5);
      \draw [line width=1pt] (-4.5,4.5)-- (-4.5,-1.5);
      \draw [line width=1pt] (-4.5,-1.5)-- (-3,-1.5);
      \draw (-3,-1.5)-- (-3,4.5);
      \draw [line width=1pt] (-3,4.5)-- (-4.5,4.5);
    \end{tikzpicture}
  \begin{tabular}{|n{2}{0}|n{2}{3} n{3}{2}|n{2}{3} n{2}{2}|n{2}{3} n{2}{2}|}
  \hline
  &
  \multicolumn{2}{c|}{$32$ bits} &
  \multicolumn{2}{c|}{$64$ bits} &
  \multicolumn{2}{c|}{$128$ bits}
  \\
  &
  \multicolumn{1}{c}{$T_c$} &
  \multicolumn{1}{c|}{$1000\times R_c$} &
  \multicolumn{1}{c}{$T_c$} &
  \multicolumn{1}{c|}{$1000\times R_c$} &
  \multicolumn{1}{c}{$T_c$} &
  \multicolumn{1}{c|}{$1000\times R_c$}
  \\
  \hline
  2 &	0.028 &	    10.22 &	 0.048 &	 15.80 &	 0.086 &	49.56 \\
  3 &	0.044 &	     6.25 &	 0.066 &	 10.77 &	 0.151 &	34.19 \\
  4 &	0.069 &	     5.18 &	 0.099 &	  9.71 &	 0.223 &	21.89 \\
  5 &	0.140 &	     3.78 &	 0.151 &	  7.30 &	 0.349 &	19.89 \\
  6 &	0.222 &	     2.91 &	 0.254 &	  6.07 &	 0.688 &	17.25 \\
  7 &	0.252 &	     2.82 &	 0.465 &	  5.06 &	 2.322 &	15.72 \\
  8 &	0.484 &	     2.44 &	 0.721 &	  4.73 &	 3.734 &	14.52 \\
  9 &	0.772 &	     2.40 &	 1.709 &	  4.37 &	 7.446 &	16.12 \\
 10 &	4.227 &	     1.83 &	 3.465 &	  4.98 &	11.364 &	16.87 \\
  \hline
  \multicolumn{7}{c}{Normalized cost}
  \\
  \hline
  &
  \multicolumn{1}{c}{$16 \times T_c$} &
  \multicolumn{1}{c|}{$16000\times R_c$} &
  \multicolumn{1}{c}{$4 \times T_c$} &
  \multicolumn{1}{c|}{$4000\times R_c$} &
  \multicolumn{1}{c}{$T_c$} &
  \multicolumn{1}{c|}{$1000\times R_c$}
  \\
  \hline
  2 &   0.444 &	  164.80 &	 0.192 &    62.90 &	 0.086 &	49.40 \\
  3 &	  0.700 &	  101.20 &	 0.264 &	  43.10 &	 0.151 &	34.06 \\
  4 &   1.112 &	   83.00 &	 0.398 &	  38.82 &	 0.223 &	23.08 \\
  5 &	  2.248 &	   60.36 &	 0.606 &    29.20 &	 0.349 &	19.89 \\
  6 &	  3.560 &	   45.80 &	 1.016 &	  24.26 &	 0.688 &	17.28 \\
  7 &	  4.028 &	   46.64 &	 1.862 &	  20.30 &	 2.322 &	15.73 \\
  8 &	  7.748 &	   38.20 &	 2.886 &	  19.20 &	 3.734 &	14.49 \\
  9 &  12.348 &	   37.64 &	 6.838 &	  17.20 &	 7.446 &	16.16 \\
 10 &  67.640 &	   29.32 &	13.862 &	  20.04 &	11.364 &	16.83 \\
  \hline
  \multicolumn{7}{|c|}{Time measured in microseconds}
  \\
  \hline
  \end{tabular}
  \end{center}
\end{table}
\ \\
Once the size $b$ has been chosen, it should be more convenient to choose the size $c$ 
of the tables used for the M4RM algorithm.
Since the cost to perform the product of the two matrices 
$\Amat \in \matSet(n,b,\FF_2)$ and $\Bmat \in \matSet(b,b,\FF_2)$
is given by \eqref{eq:M4RM:cost}, the optimal table size $c$ (when $b$ and $n$ are known)
can be estimated minimizing \eqref{eq:M4RM:cost}, which is the minimum of the following function
\begin{EQ}
  C(b,c,n) = \left(\frac{b}{c}\right)(2^{c}-1+n).
\end{EQ}
\begin{remark}
Actually, to determine the minimum of the function $C$ is quite complicated. 
However, it can be observed that $C$, as a function of $n$, 
is a straight line with a slope that decreases when $c$ increases.
So, when $n$ exceeds the value $C (b, c, n) = C (b, c +1, n)$, 
it is convenient to use a table of size $c+1$ instead of size $c$. 
In the case of the product of square matrices, the minimum cost is obtained
minimizing $C (b, c, b)$ (with respect to $c$) and we obtain the following values:
\begin{EQ}
   \arg\min_{c}\{C(32,c,32)\} \approx 4.08 \qquad
   \arg\min_{c}\{C(64,c,64)\} \approx  4.77 \\
   \arg\min_{c}\{C(128,c,128)\} \approx  5.5 \,.
\end{EQ}
\end{remark}

In Table~\ref{tab1:rows and tables costs} the costs to perform rows operations and table constructions are given;
in Table~\ref{tab2:costoprodotto} we report the costs to perform product of square matrices of size $32$, $64$ and $128$.  

\begin{table}[!htcb]
  \caption{Cost of the product of $2$ square matrices having size respectively $32$, $64$, $128$. The normalized costs are also reported.}
  \label{tab2:costoprodotto}
  \begin{center}
    \begin{tikzpicture}[scale=0.4,line cap=round,line join=round,>=triangle 45,x=1.0cm,y=1.0cm]
      \clip(-4.3,2.4) rectangle (3.5,6.3);
      \fill[fill=black,pattern=horizontal lines] (1.5,4.5) -- (1.5,3) -- (3,3) -- (3,4.5) -- cycle;
      \fill[fill=black,pattern=horizontal lines] (-4,4.5) -- (-4,3) -- (-2.5,3) -- (-2.5,4.5) -- cycle;
      \fill[fill=black,pattern=horizontal lines] (-0.5,4.5) -- (-0.5,3) -- (1,3) -- (1,4.5) -- cycle;
      \draw [line width=1pt] (1.5,4.5)-- (1.5,3);
      \draw [line width=1pt] (1.5,3)-- (3,3);
      \draw [line width=1pt] (3,3)-- (3,4.5);
      \draw [line width=1pt] (3,4.5)-- (1.5,4.5);
      \draw (-0.5,6.32) node[anchor=north west] {\Large$\bm{A}$};
      \draw (1.5,6.28) node[anchor=north west]  {\Large$\bm{B}$};
      \draw (-3.8,6.18) node[anchor=north west] {\Large$\bm{C}$};
      \draw (-2.2,4.18) node[anchor=north west] {\Large$\bm{=}$};
      \draw [line width=1pt] (-4,4.5)-- (-4,3);
      \draw [line width=1pt] (-4,3)-- (-2.5,3);
      \draw [line width=1pt] (-2.5,3)-- (-2.5,4.5);
      \draw [line width=1pt] (-2.5,4.5)-- (-4,4.5);
      \draw [line width=1pt] (-0.5,4.5)-- (-0.5,3);
      \draw [line width=1pt] (-0.5,3)-- (1,3);
      \draw [line width=1pt] (1,3)-- (1,4.5);
      \draw [line width=1pt] (1,4.5)-- (-0.5,4.5);
    \end{tikzpicture}
    \qquad\qquad\qquad
    \begin{tikzpicture}[scale=0.4, line cap=round,line join=round,>=triangle 45,x=1.0cm,y=1.0cm]
      \clip(-4.74,0.8) rectangle (6.86,6.5);
      \fill[fill=black,pattern=horizontal lines] (3.5,4.5) -- (3.5,3) -- (5,3) -- (5,4.5) -- cycle;
      \fill[fill=black,pattern=horizontal lines] (5,4.5) -- (5,3) -- (6.5,3) -- (6.5,4.5) -- cycle;
      \fill[fill=black,pattern=horizontal lines] (3.5,3) -- (3.5,1.5) -- (5,1.5) -- (5,3) -- cycle;
      \fill[fill=black,pattern=horizontal lines] (5,3) -- (5,1.5) -- (6.5,1.5) -- (6.5,3) -- cycle;
      \fill[fill=black,pattern=horizontal lines] (-4.5,4.5) -- (-4.5,3) -- (-3,3) -- (-3,4.5) -- cycle;
      \fill[fill=black,pattern=horizontal lines] (-3,4.5) -- (-3,3) -- (-1.5,3) -- (-1.5,4.5) -- cycle;
      \fill[fill=black,pattern=horizontal lines] (-4.5,3) -- (-4.5,1.5) -- (-3,1.5) -- (-3,3) -- cycle;
      \fill[fill=black,pattern=horizontal lines] (-3,3) -- (-3,1.5) -- (-1.5,1.5) -- (-1.5,3) -- cycle;
      \fill[fill=black,pattern=horizontal lines] (0,4.5) -- (0,3) -- (1.5,3) -- (1.5,4.5) -- cycle;
      \fill[fill=black,pattern=horizontal lines] (1.5,4.5) -- (1.5,3) -- (3,3) -- (3,4.5) -- cycle;
      \fill[fill=black,pattern=horizontal lines] (0,3) -- (0,1.5) -- (1.5,1.5) -- (1.5,3) -- cycle;
      \fill[fill=black,pattern=horizontal lines] (1.5,3) -- (1.5,1.5) -- (3,1.5) -- (3,3) -- cycle;
      \draw [line width=1pt] (3.5,4.5)-- (3.5,3);
      \draw [line width=1pt] (3.5,3)-- (5,3);
      \draw [line width=1pt] (5,3)-- (5,4.5);
      \draw [line width=1pt] (5,4.5)-- (3.5,4.5);
      \draw (0.8,6.18) node[anchor=north west] {\Large$\bm{A}$};
      \draw (4,6.18) node[anchor=north west] {\Large$\bm{B}$};
      \draw (-3.8,6.18) node[anchor=north west] {\Large$\bm{C}$};
      \draw (-1.44,3.5) node[anchor=north west] {\Large$\bm{=}$};
      \draw [line width=1pt] (5,4.5)-- (5,3);
      \draw [line width=1pt] (5,3)-- (6.5,3);
      \draw [line width=1pt] (6.5,3)-- (6.5,4.5);
      \draw [line width=1pt] (6.5,4.5)-- (5,4.5);
      \draw [line width=1pt] (3.5,3)-- (3.5,1.5);
      \draw [line width=1pt] (3.5,1.5)-- (5,1.5);
      \draw [line width=1pt] (5,1.5)-- (5,3);
      \draw (5,3)-- (3.5,3);
      \draw [line width=1pt] (5,3)-- (5,1.5);
      \draw [line width=1pt] (5,1.5)-- (6.5,1.5);
      \draw [line width=1pt] (6.5,1.5)-- (6.5,3);
      \draw (6.5,3)-- (5,3);
      \draw [line width=1pt] (-4.5,4.5)-- (-4.5,3);
      \draw [line width=1pt] (-4.5,3)-- (-3,3);
      \draw [line width=1pt] (-3,3)-- (-3,4.5);
      \draw [line width=1pt] (-3,4.5)-- (-4.5,4.5);
      \draw [line width=1pt] (-3,4.5)-- (-3,3);
      \draw [line width=1pt] (-3,3)-- (-1.5,3);
      \draw [line width=1pt] (-1.5,3)-- (-1.5,4.5);
      \draw [line width=1pt] (-1.5,4.5)-- (-3,4.5);
      \draw [line width=1pt] (-4.5,3)-- (-4.5,1.5);
      \draw [line width=1pt] (-4.5,1.5)-- (-3,1.5);
      \draw [line width=1pt] (-3,1.5)-- (-3,3);
      \draw (-3,3)-- (-4.5,3);
      \draw [line width=1pt] (-3,3)-- (-3,1.5);
      \draw [line width=1pt] (-3,1.5)-- (-1.5,1.5);
      \draw [line width=1pt] (-1.5,1.5)-- (-1.5,3);
      \draw (-1.5,3)-- (-3,3);
      \draw [line width=1pt] (0,4.5)-- (0,3);
      \draw [line width=1pt] (0,3)-- (1.5,3);
      \draw [line width=1pt] (1.5,3)-- (1.5,4.5);
      \draw [line width=1pt] (1.5,4.5)-- (0,4.5);
      \draw [line width=1pt] (1.5,4.5)-- (1.5,3);
      \draw [line width=1pt] (1.5,3)-- (3,3);
      \draw [line width=1pt] (3,3)-- (3,4.5);
      \draw [line width=1pt] (3,4.5)-- (1.5,4.5);
      \draw [line width=1pt] (0,3)-- (0,1.5);
      \draw [line width=1pt] (0,1.5)-- (1.5,1.5);
      \draw [line width=1pt] (1.5,1.5)-- (1.5,3);
      \draw (1.5,3)-- (0,3);
      \draw [line width=1pt] (1.5,3)-- (1.5,1.5);
      \draw [line width=1pt] (1.5,1.5)-- (3,1.5);
      \draw [line width=1pt] (3,1.5)-- (3,3);
      \draw (3,3)-- (1.5,3);
    \end{tikzpicture}
  \begin{tabular}{|n{2}{0}|n{2}{4}n{5}{4}|n{2}{4}n{5}{4}|n{3}{4}|}
  \hline 
  & 
  \multicolumn{2}{c|}{$32$ bit} &
  \multicolumn{2}{c|}{$64$ bit} &
  \multicolumn{1}{c|}{$128$ bit}
  \\
  \multicolumn{1}{|c|}{bit} & 
  \multicolumn{1}{c}{$32\times 32$} &
  \multicolumn{1}{c|}{normalized $(\times 64)$} &
  \multicolumn{1}{c}{$64\times 64$} &
  \multicolumn{1}{c|}{normalized $(\times 8)$} &
  \multicolumn{1}{c|}{$128\times 128$} \\
  \hline
  2 & 0.362 &	 22.88  &  1.060 &	 8.48  &	  6.43 \\
  3 & 0.255 &	 16.16  &  0.755 &	 6.04  &	  4.51 \\
  4 & 0.240 &	 15.20* &  0.725 &	 5.80  &	  2.88 \\
  5 & 0.262 &	 16.64  &  0.625 &	 5.00* &    2.90 \\
  6 & 0.312 &	 20.00  &  0.645 &	 5.12  &	  2.89* \\
  7 & 0.345 &	 22.08  &  0.790 &	 6.32  &	  3.91 \\
  8 & 0.570 &	 36.64  &  1.025 &	 8.16  &	  5.52 \\
  9 & 0.845 &	 54.40  &  1.955 &   15.60 &	  9.57 \\
 10 & 4.355 &	279.80  &  3.635 &	 29.08 &	 13.56 \\
  \hline
  \end{tabular}
  \end{center}
\end{table}

\subsection{How to choose switching point for Strassen Matrix-Matrix multiplication}
\label{sec:switching point}
Starting by the formula \eqref{eq:M4RM:cost}, we can easily obtain the cost to multiply (using the M4RM algorithm) two matrices of size $n$.
Let $N=\left\lceil n/b\right\rceil$ be the number of blocks in which we divide the matrix. We have to apply $N^2$ times the M4RM algorithm and we obtain the following cost  
\begin{EQ}
  \texttt{M}(n) = N^2(\texttt{T}_c^{b} + n \texttt{R}_c^{b})
  =\frac{\Xor_{b}}{b\,c}\left(n^2(2^{c}-1)+n^3\right)\, .
\end{EQ}
To perform Strassen matrix-matrix multiplication algorithm, our implementation
needs essentially $22$ additions and $7$ multiplications of matrices having size $\frac{n}{2}$.
Since the cost to compute one addition is given by $\left(\frac{n N}{4}\right)\Xor_{b}$, the whole cost for Strassen's algorithm is
\begin{EQ}
  \texttt{S}(n) = \frac{11 n^2}{2b} \Xor_{b}+ 7\,\texttt{S}\left(\frac{n}{2}\right)\, .
\end{EQ}
Then, it is convenient to use Strassen's algorithm as long as $\texttt{S}(n)\leq \texttt{M}(n)$ or
\begin{EQ}
  \texttt{S}(n) 
  = \frac{11n^2}{2b} \Xor_{b}+7\,\texttt{S}(n/2)
  \leq
  \frac{11n^2}{2b} \Xor_{b}+7\,\texttt{M}(n/2)
  \leq \texttt{M}(n) \, .
\end{EQ}
In other words, it is convenient to use Strassen's algorithm when 
\begin{EQ}
  n \geq 44 c + 6(2^{c}-1)\, .
\end{EQ}

\subsection{Complexity estimation}
\label{sec:complexity estimation}
Let $\Amat \in \matSet(n,m,\FF_2)$ be a rectangular dense matrix with $n=b N$ and $m=bM$. 
Our matrix decomposition exclusively performs row operations acting on 
string of $b$-bits.
In order to give an estimation of the complexity of our algorithm, 
we have to consider the cost to build the pseudo-inverse $\Zmat$ 
and the cost to perform our implementation of the M4RM algorithm.
For simplicity, we assume that the cost to perform one $\Xor_b$, one $\And_b$ 
and one assign operation is always the same. 
We denote by $\popC$ the cost to perform the population count operation. \\
We consider only the case $n\geq m$
to simplify exposition, even if the case $n<m$
is quite similar. Moreover the algorithm is faster when applied to matrices
with more rows than columns, so that in this case it is convenient
to apply the algorithm to the transposed matrix.
According to the procedure to build the pseudo-inverse $\Zmat$ described 
in Section~\ref{sec:incremental}, the costs, in term of operations acting on $b$ bits,
are resumed in Table~\ref{tab:CostoBuildZ}.
     
\begin{table}[!htcb]
  \caption{Cost of the BuildZ procedure considering blocks of $b$ bits}
  \label{tab:CostoBuildZ}
    \begin{center}
    \begin{tabular}{|c|c|}
    \hline 
    Line & Cost \\
    \hline
    2 & $2b$ \\
    4 & $b$ \\
    6 & $b(b-1)$  \\
    9 & $5b + rb \, {\tt popC}$ \\
    13 & $b$ \\
    15 & $3b(b-1)/2$ \\
    18 & $3b(b-1)/2$ \\
    24 & $4b^2$\\
    \hline
  \end{tabular}
  \end{center}
\end{table}
Using Table~\ref{tab:CostoBuildZ} the total cost to build all the matrices $\Zmat$ 
along the whole decomposition is given by the following formula 
\begin{EQ}\label{eq:cost:Z}
    M(5b+8b^2)+b \popC \sum_{k=1}^{M}{Nb}
    = m(5+8b)+ mn \popC = \mathcal{O}(n^2).
\end{EQ}
As seen in Section \ref{sec:6}, the cost to perform the product, using M4RM algorithm, of a matrix in $\matSet(n,b,\FF_2)$ by a matrix in $\matSet(b,b,\FF_2)$ is given by (\ref{eq:M4RM:cost}). The cost to perform the M4RM algorithm to a matrix in $\matSet(n,b,\FF_2)$ by a matrix in $\matSet(b,m,\FF_2)$ becomes 
\begin{EQ}
  \left(\frac{m}{c}\right)(2^{c}-1+n),
\end{EQ}
where the parameter $c$ is the size of the used tables.
Noticing that $m=Mb$ and $n=Nb$, we analyze the cost in term of operations
on $b$ bits in the following two extreme situations.
\begin{itemize}
  \item[a)] When $\Amat$ is full rank the cost for the M4RM algorithm applied into 
  the factorization:
  \begin{EQ}\label{M4RM:MaxRank}
    \textrm{cost}_{\mathrm{full}} =
    \sum_{k=1}^{M} \left(\frac{b}{c}\right)(M-k+1)\left(2^c-1+(N-k)b\right)
    = \frac{n^3}{3bc}+\mathcal{O}(n^2).
  \end{EQ}
  \item[b)] When $\rk(\Amat)=1$ the cost is given by
    \begin{EQ}[rcl]\label{M4RM:zeroRank}
      \textrm{cost}_{1} &=& 
      \sum_{k=1}^{M} \left(\frac{b}{c}\right)(M-k+1)(2^c-1+ Nb)
       = \frac{n^3}{2bc}+\mathcal{O}(n^2).
    \end{EQ}
\end{itemize}

This expansion holds when the table size $c$ and $b$ are fixed.
In case of $c\equiv c(n)$ or when $b\equiv b(n)$
estimation holds if the growth is asymptotically bounded by
$\lim_{n\to\infty}c(n)/\log n < \infty$ and $\lim_{n\to\infty}b(n)/n=0$.
Thus, from~\eqref{eq:cost:Z} the asymptotic cost of the factorization
is dominated by the leading term $n^3/(3bc)$ for the full rank case
and $n^3/(2bc)$ when $\rk(\Amat)=1$.
This is the cost in term of number of operations on blocks of $b$ bits. 
Andr\'en, Hellstrom and Markstrom in \cite{Markstrom2007} analyzed the cost of 
factorization in terms of number of rows operation.
Thus, in order to compare our cost with the one in this reference, 
we have to multiply our cost by $b/n$ obtaining
\begin{EQ}
  \textrm{row-cost}_{\mathrm{full}} = \frac{n^2}{3c}+\mathcal{O}(n),\qquad
  \textrm{row-cost}_{1} = \frac{n^2}{2c}+\mathcal{O}(n).
\end{EQ}
Assuming as in~\cite{Markstrom2007} the table size $c=\log_2 n$, we obtain 
\begin{EQ}
  \textrm{row-cost}_{\mathrm{full}} = \frac{n^2}{3\log_{2}{n}}+\mathcal{O}(n),\qquad
  \textrm{row-cost}_{1}=\frac{n^2}{2\log_{2}{n}}+\mathcal{O}(n).
\end{EQ}
Notice that the theoretical minimum for the complexity obtained in \cite{Markstrom2007}  
is $n^2/(2\log_{2}{n})+\mathcal{O}(n)$ and it apparently contradicts estimation
$\textrm{row-cost}_{\mathrm{full}}$.
But this minimum is obtained considering row operations performed on the whole row
while we use operations on strings of $b$-bits, so that there is no contradiction.
Our algorithm in the worst case is at least asymptotically twice faster
than the one described in~\cite{Markstrom2007}, although our algorithm 
needs to store the whole M4RM tables while Andr\'en et al. does not need additional 
memory.

\begin{remark}
In practice, the worst case is never reached because we have adopted an implementation strategy
which gives us an advantage in the case of matrices having very low rank, and in particular when the rank is equal to one.
\end{remark}


\section{Performance tests and comparison}\label{sec:7}
To evaluate the performance of our algorithm, we compute the decomposition of $n \times n$ matrices with entries in $\FF_2$.
Our block decomposition works well both for dense matrices and for relatively sparse matrices.
Notice that in the former case the obtained matrices have rank most probably equal to $n-1$ (or $n$);
in the latter case we have low-rank matrices. 
  
First, we have constructed a sample of random dense matrices of size $n$ 
(where the size $n$ ranges from $256$ and $65536$).
   
In Table~\ref{OurTimesComparison} we give the minimum value of ten observed running times for computing our block decomposition (recursive and non-recursive), in case $b=32, 64, 128$.

\begin{table}[!htcb]
  \caption{Compare running times with $b=32$, $64$, $128$ recursive (rec) and non 
           recursive versions of our algorithm.
           Execution obtained using LLVM 3.0 compiler (MAC OSX) on a 3.06GHz Intel(R) Xeon.}
  \begin{center}
  \begin{tabular}{|c||n{3}{3}|n{3}{3}|n{3}{3}|n{3}{3}|n{3}{3}|n{3}{3}|}
  \hline
  \multicolumn{1}{|c||}{}                & 
  \multicolumn{1}{c|}{$32$}              & 
  \multicolumn{1}{c|}{\textrm{$32$ rec}} & 
  \multicolumn{1}{c|}{$64$}              &
  \multicolumn{1}{c|}{\textrm{$64$ rec}} &
  \multicolumn{1}{c|}{$128$}             &
  \multicolumn{1}{c|}{\textrm{$128$ rec}} \\
  \hline
  \multicolumn{1}{|c||}{$n$} & \multicolumn{6}{c|}{\textrm{milliseconds}}\\
  \hline
  256  &  0.106 &  0.108 &  0.131 &  0.133 &  0.247 &  0.250 \\ 
  512  &  0.347 &  0.350 &  0.327 &  0.330 &  0.553 &  0.557 \\
  1024 &  1.574 &  1.741 &  1.070 &  1.171 &  1.443 &  1.527 \\
  2048 &  9.260 & 10.862 &  5.076 &  5.927 &  5.172 &  5.795 \\
  4096 & 64.849 & 76.458 & 32.454 & 37.996 & 27.274 & 30.792 \\
  \hline
  \multicolumn{1}{|c||}{$n$} & \multicolumn{6}{c|}{\textrm{seconds}}\\
  \hline
   8192 &  0.492 &  0.549 &  0.248 &  0.267 &  0.203 &  0.206\\ 
  16384 &  3.922 &  3.981 &  1.939 &  1.981 &  1.646 &  1.523\\
  32768 & 32.011 & 29.791 & 23.859 & 15.440 & 12.884 & 11.498\\
  65536 & 253.169 & 216.323 & 190.065 & 114.794 & 102.071 & 86.156\\
  \hline
  \end{tabular}
  \end{center}
  \label{OurTimesComparison}
\end{table}

In Table \ref{CompareRunninTimes}, the minimum running time of ten trials to obtain a matrix decomposition is given. In particular, we compare the running times to obtain M4RI, PLUQ, PLE  decompositions adopted into SAGE~\cite{CGC-misc-sage} with the ones to get our $128$ block recursive decomposition.

\begin{table}[!htcb]
  \caption{Compare running times with three algorithms (PLUQ, M4RI, PLE) adopted into SAGE
           and our recursive algorithm with $b=128$ 
           using LLVM 3.0 compiler (MAC OSX) on a 3.06GHz Intel(R) Xeon.}
  \begin{center}
  \begin{tabular}{|n{5}{0}||n{3}{3}|n{3}{3}|n{3}{3}|n{3}{3}|}
  \hline
     & \hfil\textrm{M4RI} & \hfil\textrm{PLUQ} & \hfil\textrm{PLE} & \hfil\textrm{our}\\
  \hline
  \hfil $n$ & \multicolumn{4}{c|}{\textrm{milliseconds}}\\
  \hline
  256   &   0.248  &   0.239 &   0.239 &  0.252 \\
  512   &   0.847  &   0.836 &   0.821 &  0.562 \\
  1024  &   3.428  &   2.653 &   2.551 &  1.539 \\
  2048  &  12.810  &  10.509 &  10.003 &  5.850 \\
  4096  &  62.231  &  53.686 &  50.400 & 31.406 \\
  \hline
  \hfil $n$ & \multicolumn{4}{c|}{\textrm{seconds}}\\
  \hline
  8192  &   0.349 &   0.340 &   0.334 &  0.207 \\ 
  16384 &   2.462 &   2.425 &   2.392 &  1.532 \\
  32768 &  18.674 &  18.558 &  18.341 & 11.529 \\
  65536 & 135.906 & 128.567 & 125.783 & 86.538 \\
  \hline
  \end{tabular}
  \end{center}
    \label{CompareRunninTimes}
\end{table}
Then, we constructed low-rank matrices as follows. 
We consider the samples 
$$S_n=\{S_{n, i} \,| \, i =2, \ldots, n\}$$ 
of $39$ relatively sparse matrices of size $n$ having respectively $i$ non-zero elements per rows. 

\begin{figure}[!htcb]
  \begin{center}
    \begin{tabular}{rr}
    \pgfplotstableread{nnz1024.dat}{\datafile}
    \begin{tikzpicture}[scale=1.1]
    \begin{axis}[ width=6cm,
                  height=4cm,
                  legend columns=4, legend style={draw=none},
                  minor tick num=1,
                  title=$1024\times1024$,
                  xmin = 0, xmax = 40]
    \addplot [only marks, mark=triangle*, mark size = 1.2, mark options={fill=black}] table [x={nnz}, y={mio}]  {\datafile};
    \addplot [only marks, mark=o,         mark size = 1.2 ] table [x={nnz}, y={PLE}]  {\datafile};
    \addplot [only marks, mark=square,    mark size = 1.2 ] table [x={nnz}, y={PLU}]  {\datafile};
    \addplot [only marks, mark=+,         mark size = 1.2 ] table [x={nnz}, y={M4RI}] {\datafile};
    \legend{our\\PLE\\PLU\\M4RI\\}
    \end{axis}
    \end{tikzpicture}
    &
    \pgfplotstableread{nnz2048.dat}{\datafile}
    \begin{tikzpicture}[scale=1.1]
    \begin{axis}[ width=6cm,
                  height=4cm,
                  legend columns=4, legend style={draw=none},
                  minor tick num=1,
                  title=$2048\times2048$,
                  xmin = 0, xmax = 40]
    \addplot [only marks, mark=triangle*, mark size = 1.2, mark options={fill=black}] table [x={nnz}, y={mio}]  {\datafile};
    \addplot [only marks, mark=o,         mark size = 1.2 ] table [x={nnz}, y={PLE}]  {\datafile};
    \addplot [only marks, mark=square,    mark size = 1.2 ] table [x={nnz}, y={PLU}]  {\datafile};
    \addplot [only marks, mark=+,         mark size = 1.2 ] table [x={nnz}, y={M4RI}] {\datafile};
    \legend{our\\PLE\\PLU\\M4RI\\}
    \end{axis}
    \end{tikzpicture}
    \\[1em]
    \pgfplotstableread{nnz4096.dat}{\datafile}
    \begin{tikzpicture}[scale=1.1]
    \begin{axis}[ width=6cm,
                  height=4cm,
                  legend columns=4, legend style={draw=none},
                  minor tick num=1,
                  title=$4096\times4096$,
                  xmin = 0, xmax = 40]
    \addplot [only marks, mark=triangle*, mark size = 1.2, mark options={fill=black}] table [x={nnz}, y={mio}]  {\datafile};
    \addplot [only marks, mark=o,         mark size = 1.2 ] table [x={nnz}, y={PLE}]  {\datafile};
    \addplot [only marks, mark=square,    mark size = 1.2 ] table [x={nnz}, y={PLU}]  {\datafile};
    \addplot [only marks, mark=+,         mark size = 1.2 ] table [x={nnz}, y={M4RI}] {\datafile};
    \legend{our\\PLE\\PLU\\M4RI\\}
    \end{axis}
    \end{tikzpicture}
    &
    \pgfplotstableread{nnz8192.dat}{\datafile}
    \begin{tikzpicture}[scale=1.1]
    \begin{axis}[ width=6cm,
                  height=4cm,
                  legend columns=4, legend style={draw=none},
                  minor tick num=1,
                  title=$8192\times8192$,
                  xmin = 0, xmax = 40]
    \addplot [only marks, mark=triangle*, mark size = 1.2, mark options={fill=black}] table [x={nnz}, y={mio}]  {\datafile};
    \addplot [only marks, mark=o,         mark size = 1.2 ] table [x={nnz}, y={PLE}]  {\datafile};
    \addplot [only marks, mark=square,    mark size = 1.2 ] table [x={nnz}, y={PLU}]  {\datafile};
    \addplot [only marks, mark=+,         mark size = 1.2 ] table [x={nnz}, y={M4RI}] {\datafile};
    \legend{our\\PLE\\PLU\\M4RI\\}
    \end{axis}
    \end{tikzpicture}
    \\[1em]
    \pgfplotstableread{nnz16384.dat}{\datafile}
    \begin{tikzpicture}[scale=1.1]
    \begin{axis}[ width=6cm,
                  height=4cm,
                  legend columns=4, legend style={draw=none},
                  minor tick num=1,
                  title=$16384\times16384$,
                  xmin = 0, xmax = 40]
    \addplot [only marks, mark=triangle*, mark size = 1.2, mark options={fill=black}] table [x={nnz}, y={mio}]  {\datafile};
    \addplot [only marks, mark=o,         mark size = 1.2 ] table [x={nnz}, y={PLE}]  {\datafile};
    \addplot [only marks, mark=square,    mark size = 1.2 ] table [x={nnz}, y={PLU}]  {\datafile};
    \addplot [only marks, mark=+,         mark size = 1.2 ] table [x={nnz}, y={M4RI}] {\datafile};
    \legend{our\\PLE\\PLU\\M4RI\\}
    \end{axis}
    \end{tikzpicture}
    &
    \pgfplotstableread{nnz32768.dat}{\datafile}
    \begin{tikzpicture}[scale=1.1]
    \begin{axis}[ width=6cm,
                  height=4cm,
                  legend columns=4, legend style={draw=none},
                  minor tick num=1,
                  title=$32768\times32768$,
                  xmin = 0, xmax = 40]
    \addplot [only marks, mark=triangle*, mark size = 1.2, mark options={fill=black}] table [x={nnz}, y={mio}]  {\datafile};
    \addplot [only marks, mark=o,         mark size = 1.2 ] table [x={nnz}, y={PLE}]  {\datafile};
    \addplot [only marks, mark=square,    mark size = 1.2 ] table [x={nnz}, y={PLU}]  {\datafile};
    \addplot [only marks, mark=+,         mark size = 1.2 ] table [x={nnz}, y={M4RI}] {\datafile};
    \legend{our\\PLE\\PLU\\M4RI\\}
    \end{axis}
    \end{tikzpicture}
    \\[1em]
    \pgfplotstableread{nnz65536.dat}{\datafile}
    \begin{tikzpicture}[scale=1.1]
    \begin{axis}[ width=6cm,
                  height=4cm,
                  legend columns=4, legend style={draw=none},
                  minor tick num=1,
                  title=$65536\times65536$,
                  xmin = 0, xmax = 40]
    \addplot [only marks, mark=triangle*, mark size = 1.2, mark options={fill=black}] table [x={nnz}, y={mio}]  {\datafile};
    \addplot [only marks, mark=o,         mark size = 1.2 ] table [x={nnz}, y={PLE}]  {\datafile};
    \addplot [only marks, mark=square,    mark size = 1.2 ] table [x={nnz}, y={PLU}]  {\datafile};
    \addplot [only marks, mark=+,         mark size = 1.2 ] table [x={nnz}, y={M4RI}] {\datafile};
    \legend{our\\PLE\\PLU\\M4RI\\}
    \end{axis}
    \end{tikzpicture}
    \end{tabular}
  \end{center}
  \caption{Running times for computing a matrix decomposition of $n \times n$ low-rank matrices
           using LLVM 3.0 compiler (MAC OSX) on a 2.53GHz Intel(R) Core 2 Duo.}
  \label{fig:comparison}
\end{figure}
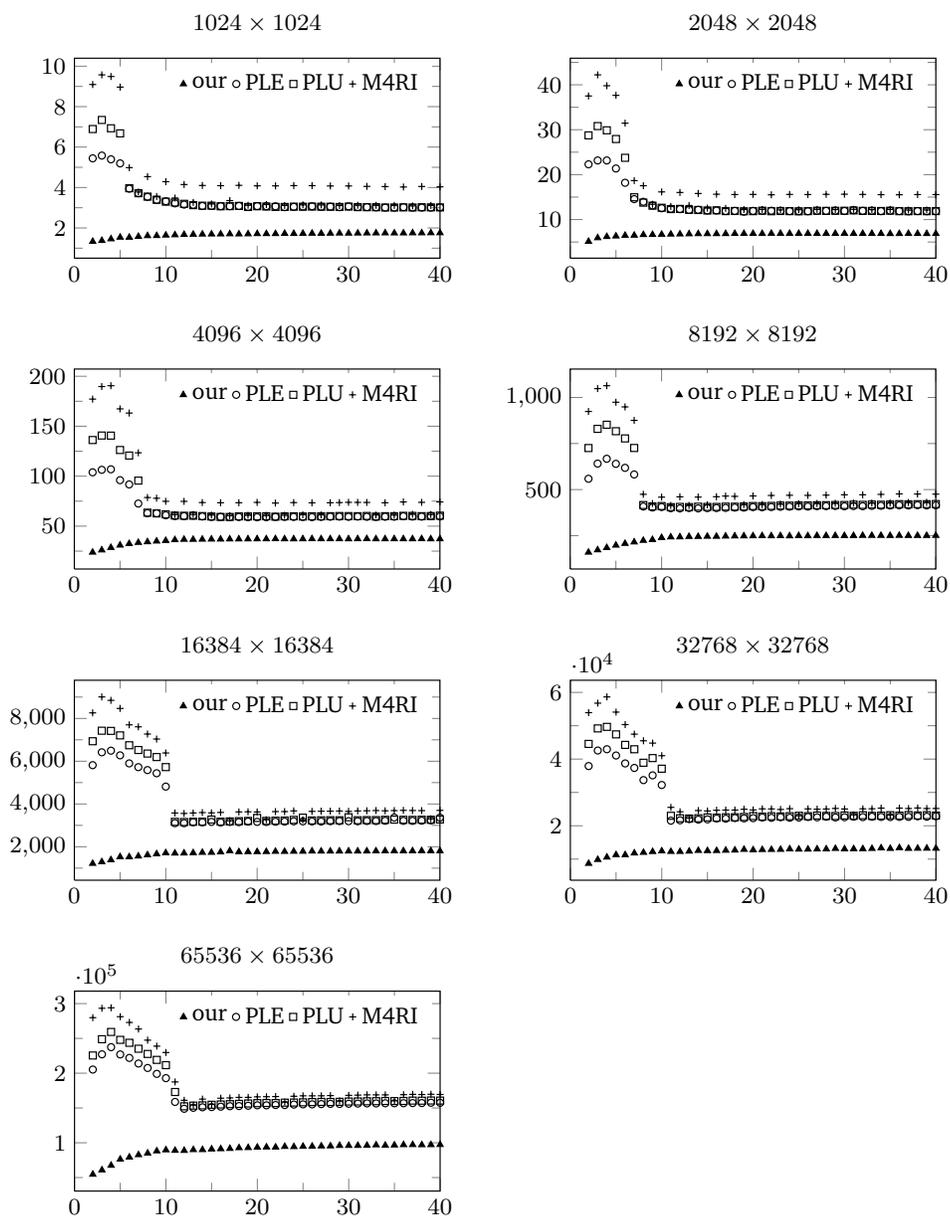

In Figure \ref{fig:comparison} we have plotted the observed running times (in milliseconds) for M4RI (cross), PLUQ (square), PLE(circles) and our $128$ bit block recursive (triangles) decomposition. This is done for every matrix in the sample $S_n$
The size $n$ ranges from $1024$ to $65536$.


\section*{Acknowledgment}
The authors would like to thank the following people for their valuable comments and suggestions: Prof. M. Elia, Prof. M. Sala and Prof. T. Muntean.\\
Part of the second author's research has been funded by the ``Groupe d'\'Etudes et Recherche en Informatique des Syst\`emes Communicants S\'ecuris\'es'', Aix-Marseille Universit\'e.

\bibliographystyle{acmtrans} 
\bibliography{RefsCGC}

\begin{thebibliography}{}

\bibitem[\protect\citeauthoryear{Aho, Hopcroft, and Ullman}{Aho
  et~al\mbox{.}}{1974}]{aho:1974}
{\sc Aho, A.~V.}, {\sc Hopcroft, J.~E.}, {\sc and} {\sc Ullman, J.~D.} 1974.
\newblock {\em The Design and Analysis of Computer Algorithms}.
\newblock Addison-Wesley Publishing Company, Boston, Massachusett.

\bibitem[\protect\citeauthoryear{Albrecht, Bard, and Hart}{Albrecht
  et~al\mbox{.}}{2010}]{Albrecht:2010}
{\sc Albrecht, M.}, {\sc Bard, G.}, {\sc and} {\sc Hart, W.} 2010.
\newblock Algorithm 898: Efficient multiplication of dense matrices over
  $\mathrm{GF}(2)$.
\newblock {\em ACM Trans. Math. Softw.\/}~{\em 37}, 9:1--9:14.

\bibitem[\protect\citeauthoryear{Albrecht, Bard, and Pernet}{Albrecht
  et~al\mbox{.}}{2011}]{Albrecht:2011}
{\sc Albrecht, M.~R.}, {\sc Bard, G.~V.}, {\sc and} {\sc Pernet, C.} 2011.
\newblock Efficient dense gaussian elimination over the finite field with two
  elements.
\newblock Tech. rep., Cornell University Library.

\bibitem[\protect\citeauthoryear{Andr\'en, Hellstrom, and Markstrom}{Andr\'en
  et~al\mbox{.}}{2007}]{Markstrom2007}
{\sc Andr\'en, D.}, {\sc Hellstrom, L.}, {\sc and} {\sc Markstrom, K.} 2007.
\newblock On the complexity of matrix reduction over finite fields.
\newblock {\em Advances in Applied Mathematics\/}~{\em 39,\/}~4, 428 -- 452.

\bibitem[\protect\citeauthoryear{Arlazarov, Dinic, Kronrod, and
  Farad{\v{z}}ev}{Arlazarov et~al\mbox{.}}{1970}]{Arlazarov:1970}
{\sc Arlazarov, V.~L.}, {\sc Dinic, E.~A.}, {\sc Kronrod, M.~A.}, {\sc and}
  {\sc Farad{\v{z}}ev, I.~A.} 1970.
\newblock On economical construction of the transitive closure of a directed
  graph.
\newblock {\em Soviet Mathematics---Doklady\/}~{\em 11,\/}~5, 1209--1210.

\bibitem[\protect\citeauthoryear{Bach and Shallit}{Bach and
  Shallit}{1996}]{bach:1996}
{\sc Bach, E.} {\sc and} {\sc Shallit, J.} 1996.
\newblock {\em Algorithmic Number Theory: Efficient algorithms}.
\newblock Foundations of computing. MIT Press, Cambridge, USA.

\bibitem[\protect\citeauthoryear{Ben-Israel and Greville}{Ben-Israel and
  Greville}{2003}]{Israel:2004}
{\sc Ben-Israel, A.} {\sc and} {\sc Greville, T.~N.} 2003.
\newblock {\em {Generalized inverses. Theory and applications. 2nd ed.}}
\newblock CMS Books in Mathematics. Springer.

\bibitem[\protect\citeauthoryear{Coppersmith and Winograd}{Coppersmith and
  Winograd}{1990}]{Coppersmith:1990}
{\sc Coppersmith, D.} {\sc and} {\sc Winograd, S.} 1990.
\newblock Matrix multiplication via arithmetic progressions.
\newblock {\em J. Symbolic Comput.\/}~{\em 9,\/}~3, 251--280.

\bibitem[\protect\citeauthoryear{Douglas, Heroux, Slishman, and Smith}{Douglas
  et~al\mbox{.}}{1994}]{Douglas:1994}
{\sc Douglas, C.~C.}, {\sc Heroux, M.}, {\sc Slishman, G.}, {\sc and} {\sc
  Smith, R.~M.} 1994.
\newblock {GEMMW:} a portable level 3 blas winograd variant of strassen's
  matrix-matrix multiply algorithm.
\newblock {\em Journal of Computational Physics\/}~{\em 110,\/}~1, 1 -- 10.

\bibitem[\protect\citeauthoryear{Golub and Van~Loan}{Golub and
  Van~Loan}{1996}]{Golub:1996}
{\sc Golub, G.~H.} {\sc and} {\sc Van~Loan, C.~F.} 1996.
\newblock {\em Matrix computations\/}, Third ed.
\newblock Johns Hopkins Studies in the Mathematical Sciences. Johns Hopkins
  University Press, Baltimore, MD.

\bibitem[\protect\citeauthoryear{Haynsworth}{Haynsworth}{1968}]{haynsworth:1968}
{\sc Haynsworth, E.} 1968.
\newblock {\em On The Schur Complement}.
\newblock Defense Technical Information Center, Basel University.

\bibitem[\protect\citeauthoryear{Higham}{Higham}{1990}]{Higham:1990}
{\sc Higham, N.~J.} 1990.
\newblock Exploiting fast matrix multiplication within the level 3 blas.
\newblock {\em ACM Trans. Math. Softw.\/}~{\em 16}, 352--368.

\bibitem[\protect\citeauthoryear{Huss-Lederman, Jacobson, Tsao, Turnbull, and
  Johnson}{Huss-Lederman et~al\mbox{.}}{1996}]{Huss-Lederman:1996}
{\sc Huss-Lederman, S.}, {\sc Jacobson, E.~M.}, {\sc Tsao, A.}, {\sc Turnbull,
  T.}, {\sc and} {\sc Johnson, J.~R.} 1996.
\newblock Implementation of strassen's algorithm for matrix multiplication.
\newblock In {\em Proceedings of the 1996 ACM/IEEE conference on Supercomputing
  (CDROM)}. Supercomputing '96. IEEE Computer Society, Washington, DC, USA.

\bibitem[\protect\citeauthoryear{Jeannerod, Pernet, and Storjohann}{Jeannerod
  et~al\mbox{.}}{2011}]{Jeannerod:2011}
{\sc Jeannerod, C.-P.}, {\sc Pernet, C.}, {\sc and} {\sc Storjohann, A.} 2011.
\newblock Rank-profile revealing gaussian elimination and the cup matrix
  decomposition.
\newblock Tech. rep., Cornell University Library.

\bibitem[\protect\citeauthoryear{Kaporin}{Kaporin}{1999}]{Kaporin:1999}
{\sc Kaporin, I.} 1999.
\newblock A practical algorithm for faster matrix multiplication.
\newblock {\em Numerical Linear Algebra with Applications\/}~{\em 6,\/}~8,
  687--700.

\bibitem[\protect\citeauthoryear{Kincaid and Cheney}{Kincaid and
  Cheney}{2002}]{kincaid:2002}
{\sc Kincaid, D.} {\sc and} {\sc Cheney, E.} 2002.
\newblock {\em Numerical Analysis: Mathematics of Scientific Computing}.
\newblock Pure and applied undergraduate texts. American Mathematical Society.

\bibitem[\protect\citeauthoryear{Meyer}{Meyer}{2000}]{Meyer:2000}
{\sc Meyer, C.~D.} 2000.
\newblock {\em Matrix analysis and applied linear algebra}.
\newblock Society for Industrial and Applied Mathematics, Philadelphia, PA,
  USA.

\bibitem[\protect\citeauthoryear{Shoup}{Shoup}{2009}]{Shoup:2009}
{\sc Shoup, V.} 2009.
\newblock {\em A computational introduction to number theory and algebra}.
\newblock Cambridge University Press.

\bibitem[\protect\citeauthoryear{Stein et~al\mbox{.}}{Stein
  et~al\mbox{.}}{2012}]{CGC-misc-sage}
{\sc Stein, W.} {\sc et~al\mbox{.}} 2012.
\newblock {\em {S}age {M}athematics {S}oftware ({V}ersion 5.2)}.
\newblock The Sage~Developement Team.
\newblock { \tt http://www.sagemath.org}.

\bibitem[\protect\citeauthoryear{Strassen}{Strassen}{1969}]{Strassen:1969}
{\sc Strassen, V.} 1969.
\newblock Gaussian elimination is not optimal.
\newblock {\em Numerische Mathematik\/}~{\em 13}, 354--356.
\newblock 10.1007/BF02165411.

\bibitem[\protect\citeauthoryear{Winograd}{Winograd}{1971}]{Winograd:1971}
{\sc Winograd, S.} 1971.
\newblock On multiplication of $2\times 2$ matrices.
\newblock {\em Linear Algebra and its Applications\/}~{\em 4,\/}~4, 381 -- 388.

\bibitem[\protect\citeauthoryear{Zhang}{Zhang}{2005}]{zhang:2005}
{\sc Zhang, F.} 2005.
\newblock {\em The Schur Complement and Its Applications}.
\newblock Numerical Methods and Algorithms. Springer Science.

\end{thebibliography}

\end{document}